\documentclass[14pt]{article}
\usepackage{amsmath,amstext,amssymb,graphicx}
\usepackage{verbatim}
\usepackage{color}
\usepackage{algorithm,algorithmic}
\usepackage{epstopdf}

\usepackage{amssymb,indentfirst,xspace,bm}
\usepackage{mathtools}

\def\Ac{{\cal A}}

\def\Gc{{\cal G}}

\def\Rc{{\cal R}}
\def\Zc{{\cal Z}}
\def\norm#1{\|#1\|}

\def\vect#1{\mbox{\boldmath{$#1$}}}

\def\mR{\mathbb{R}}
\def\mC{\mathbb{C}}

\def\wG{\widehat{G}}
\def\wU{\widehat{U}}
\def\wV{\widehat{V}}
\def\wf{\widehat{f}}
\def\we{\widehat{e}}

\def\wg{\widehat{g}}

\def\wP{\widehat P}
\def\wM{\widehat M}

\def\Rho{{X}}
\def\RhoQ{{Y}}
\def\Lc{\mathcal{L}}

\newcommand{\eps}{\varepsilon}
\renewcommand{\epsilon}{\eps}
\renewcommand{\leq}{\leqslant}

\DeclarePairedDelimiter{\abs}{\lvert}{\rvert}

\def\bfrho{\mbox{\boldmath$\rho$}}
\def\bfgamma{\mbox{\boldmath$\gamma$}}

\title{Illumination strategies for intensity-only imaging}
\author{Alexei Novikov \footnotemark[1]\ \footnotemark[4] 
\and Miguel Moscoso\footnotemark[2]\ \footnotemark[5]\ 
\and George Papanicolaou\footnotemark[3]\ \footnotemark[6]\ 
}
\begin{document}
\maketitle

\renewcommand{\thefootnote}{\fnsymbol{footnote}}
\footnotetext[1]{anovikov@math.psu.edu}
\footnotetext[2]{moscoso@math.uc3m.es}
\footnotetext[3]{papanico@math.stanford.edu}
\footnotetext[4]{Department of Mathematics, Penn State University. University Park, PA 16802, USA}
\footnotetext[5]{Gregorio Mill\'{a}n Institute, Universidad Carlos III de Madrid, Madrid 28911, Spain}
\footnotetext[6]{Department of Mathematics, Stanford University, California 94305, USA}

\renewcommand{\thefootnote}{\arabic{footnote}}

\begin{abstract}

We propose a new strategy for narrow band, active array imaging of localized scatterers 
when only the intensities are recorded and measured at the array.
We consider a homogeneous medium so that wave propagation is fully coherent. 
We show that imaging with intensity-only measurements can be carried out using the
{\em time reversal operator}  
of the imaging system, which can be obtained from intensity measurements using 
an appropriate illumination strategy and the polarization identity. 
Once the time reversal operator has been obtained, 
we show that the images can be formed using its 
singular value decomposition (SVD). We use two SVD-based methods to image the scatterers.
The proposed approach is simple and efficient. It does not need prior information 
about the sought image, and guarantees exact recovery in the noise-free case. Furthermore, it is robust with respect to additive 
noise. Detailed numerical simulations illustrate the performance of the proposed imaging strategy when only the intensities are captured.
\end{abstract}


\pagestyle{myheadings}
\thispagestyle{plain}

\section{Introduction}
Imaging using intensity-only (or phaseless) measurements is challenging because much information about the sought image 
is lost in the unrecorded phases.
The problem of recovering an image from intensity-only measurements, known as the phase retrieval problem,
arises in many situations in which it is difficult, or impossible, to measure and record the phases of the signals
received at the detectors. This is the case, for example, in imaging from X-ray sources \cite{Millane90, Harrison93, Pfeiffer06},
 or from optical sources \cite{Walther63,Dainty87,Trebino93}, where one seeks to reconstruct an image from the spectral intensities. This 
 problem arises in various fields,  including  crystallography,  optical imaging, astronomy, and electron
microscopy, and the images to be formed from intensity-only measurements vary from galaxies to microscopic objects.

In this paper, we consider the problem in active array imaging when the sensors only record  the intensities of the signals.
This can be the case because less expensive sensors are used, the data need to be collected faster, or because the phases are difficult to 
 measure at the frequencies used for imaging.
For frequencies above 10 GHz or so, it is difficult at present to record the phase of the scattered signals directly.


There are at least two different approaches for imaging using intensity-only measurements. In the first
approach, the phases are retrieved from the experimental set-up
before doing the imaging. This is done, for example, in holographic based methods
where an interferometer records the interference pattern between a reference signal and the analyzed signal 
\cite{Osherovich12,Raz13}. The interferometric image depends on the phase difference between the two signals and, hence, holds the 
desired phase information. 
An experimental strategy is also proposed for diffraction tomography in~\cite{Gbur02},
which requires measurements of the signal on two planes spaced at distances
smaller than a wavelength. Such techniques are, however, hard to implement in practice.

The second approach carries out imaging directly, without previous estimation of the missing phases,
using reconstruction algorithms. A frequently used method is based on alternating projection
algorithms, proposed  by Gerschberg and Saxton (GS)~\cite{GS72}.
This method uses two intensity measurements to form the image: the magnitude of the image itself, and the magnitude
of its Fourier transform, i.e., the spectral intensity. The GS algorithm alternates between the spatial and the frequency domains, 
correcting the current iterate by imposing constrains in the spatial domain and scaling the 
Fourier coefficients in the frequency domain. Fienup~\cite{Dainty87} proposed a successful modification
of the GS algorithm,
the Hybrid-Input-Output (HIO) algorithm, which is less prone to stagnation 
and only requires  one intensity measurement, the 
spectral intensity of the image one wishes to form. The HIO algorithm is, probably, the 
algorithm used most widely at present.  However, it is a non-convex algorithm and it does not converge in general to the exact solution, even with noiseless data. 
To increase the likelihood of convergence, 
HIO often requires  image priors (finite spatial extent, real-valuedness, positivity, etc), but this additional information is not always available. 

In \cite{Fannjiang12}, the authors propose to use a phase modulator which randomly modifies the phases of the original image by a known mask. They prove that random illuminations often lead to a unique solution  and remove the stagnation problem associated to GS and HIO algorithms. In \cite{Fannjiang13}, the uniqueness result is extended to the case where only rough information about the mask's phases is assumed.
Newton-type and other gradient-based optimization methods have also been proposed. However, these methods may fail  due to the high non-linearity
of the phase retrieval  problem~\cite{Nieto-Vesperinas86}. See also~\cite{Marchesini07} for a survey and comparison of iterative projection and gradient-based algorithms.

To overcome the problems of convergence of these algorithms, and motivated by the recent developments in compressed sensing \cite{Candes05,Donoho06}, the authors in~\cite{CMP13} proposed a convex approach that is capable of solving the  
problem of imaging using only intensities. In~\cite{CMP13},  the 
non linear vector problem in phase retrieval is replaced
by a linear matrix one, which is solved 
by  using nuclear norm  minimization.
This makes the problem convex and solvable in polynomial time, and yields the unique solution in the
noise-free case. In~\cite{Candes13}, this approach is combined with the use of masks.
They show that a few simple structured illumination patterns can determine the solution uniquely using this formulation.

While this convex approach is an important advance
for intensity-only imaging problems, 
it is computationally expensive for large scale problems, for example, for images with a large number $K$ of pixels. 
This is so, because  it requires the solution of a $K\times K$ optimization problem 
with $K^2$ unknowns, instead of the original one with $K$ unknowns. 
In other words, it transforms the phase retrieval problem into one of recovering a rank-one matrix, which leads 
to very large optimization problems that are not feasible if the images are large. As a consequence,
it is desirable to have other approaches that guarantee convergence to the exact solution and, at the
same time, keep the size of the problem small so the solution can be found more efficiently. 
It is important that any such approaches be robust to noise. 

The main contribution of this paper is the introduction of a new strategy for imaging when only the intensities are recorded.
This strategy has the desired properties mentioned above:  exact recovery, 
robustness with respect to noise, and efficiency for large problems. 
We show that imaging of a small number of localized scatterers  can be  accomplished using the
{\em time reversal operator}  $\vect\wM(\omega)=\vect\wP^*(\omega) \vect\wP(\omega)$, where $\vect\wP(\omega)$ 
is the full array response matrix of the imaging system. 
We show that the {\em time reversal operator} can be obtained from the total power recorded 
at the array using an appropriate illumination strategy and the polarization identity. 
Once the  {\em time reversal operator} has been obtained, we show that 
 the  location of the scatterers can be determined using 
its singular value decomposition (SVD).
 
We consider two methods that make use of the SVD of $\vect\wM(\omega)$.
The first method finds the locations of the scatterers 
from the perspective of sparse optimization, using a Multiple Measurement Vector (MMV) approach.
The second method  finds the locations of the  scatterers by beamforming. We use the MUSIC (MUltiple
SIgnal Classification) method, which is equivalent to beamforming,  using the significant singular vectors as illuminations. Both methods
recover the location of the  scatterers exactly in the noise-free case and are  robust with respect to additive noise.

The imaging methods described here 
are efficient, do not  need prior information about the object to be imaged,
and guarantee exact recovery. We note, however, that recording all the intensities needed for the {\em time reversal operator} 
may not be possible. Indeed, the number of illuminations involved is $N^2$, where $N$ is the number of transducers in the array. In order to simplify the data acquisition process,  we also
propose two methods that reduce the number of illuminations needed for imaging. The first method selects pairs of transducers randomly, 
and finds the missing entries in the {\em time reversal operator} via matrix completion. This method reduces the number of illuminations 
to one half. The second method does not select the transducers 
randomly, but uses only a few transducers at the edges of the array. This method reduces the number of illuminations even more.

The paper is organized as follows. In Section \ref{sec:model}, we formulate the active array imaging problem using intensity-only
measurements. In Section \ref{sec:timereversalop}, we show how to obtain the
{\em time reversal operator} when only the intensities of the signals are recorded at the array, and we discuss the relation 
of the {\em time reversal operator} with the full data matrix (that also contains the information about the phases of the signals). 
We also discuss in Section \ref{sec:timereversalop} imaging with an incomplete set of illuminations, i.e., when some entries
of the {\em time reversal operator} are missing.
In Section \ref{sec:methods}, we briefly review MMV and MUSIC methods, the two imaging methods used in the paper to form the images.
In Section \ref{sec:numerics}, we show the results of numerical experiments. Section \ref{sec:conclusions} contains our conclusions.

\section{Active array imaging}
\label{sec:model}
In active array imaging we seek to locate the positions and reflectivities of a set of scatterers using the data recorded on an array 
${\cal A}$. By an active array, we mean a collection of $N$ transducers that emit spherical wave signals from positions $\vect x_s\in{\cal A}$ and record the echoes with 
receivers at positions $\vect x_r\in{\cal A}$. The transducers are placed at distance $h$ between them, which is of the order of the wavelength 
$\lambda=2\pi c_0/\omega$, where $c_0$ is the wave speed in the medium and $\omega$ is the frequency of the probing signal. In this paper, we focus on imaging of 
localized scatterers, which means that the scatterers are very small compared to the wavelength (point-like scatterers). Furthermore, for ease of exposition, we assume that multiple scattering between the scatterers is negligible. 
The imaging methods considered here can be implemented when multiple scattering is important too (see \cite{Chai14}
for details). 

Let the active array with $N$ transducers at positions $\vect x_s$, $s=1,\cdots,N$,
be located on the plane $z=0$. Assume that there are $M$ point-like scatterers in a  image window (IW), which is at a distance $L$ from the array. We discretize the IW using a uniform grid of $K \gg M$ points $\vect y_j$, $j=1,\ldots,K$.
The scatterers have reflectivities $\alpha_j\in\mC$, and are
located at  positions $\vect y_{n_1},\ldots,\vect y_{n_M}$, which we assume coincide with one of these $K$ grid points.
If the scatterers are far apart or the reflectivities are small, interaction between scatterers is
weak and multiple scattering can be neglected.
Then, with the Born approximation,
the response at $\vect x_r$ due to
a narrow-band pulse of angular frequency $\omega$
sent from $\vect x_s$ and reflected by the $M$ scatterers is given by
\begin{equation}
\label{response}
\wP(\vect x_r,\vect x_s,\omega)=\sum_{j=1}^M\alpha_j\wG_0(\vect
x_r,\vect y_{n_j},\omega)\widehat G_0(\vect y_{n_j},\vect
x_s,\omega)\, ,
\end{equation}
where 
\begin{equation}
\label{homo_green}
\wG_0(\vect x,\vect y,\omega) = \frac{\exp\{i\kappa\abs{\vect x - \vect y}\}}{4\pi \abs{\vect x - \vect y}} 
\end{equation}
is the Green's function that characterizes 
wave propagation from $\vect x$ to $\vect y$ in a homogeneous medium. 
To write the data received on the array in a more compact form, we define the {\em Green's function vector} 
$\vect\wg_0(\vect y,\omega)$ at location $\vect y$ in IW as 
\begin{equation}
\label{GreenFuncVec}
\vect\wg_0(\vect y,\omega)=[\wG_0(\vect x_{1},\vect y,\omega),\cdots,\wG_0(\vect x_{N},\vect y,\omega)]^T\, ,
\end{equation}
where $.^T$ means the transpose. This vector can also be interpreted as the illumination vector of the array targeting the position $\vect y$. 
We also define the true {\it reflectivity vector}
$\vect\rho_0=[\rho_{01},\ldots,\rho_{0K}]^T\in\mC^K$ such that
\begin{equation}
\label{eq:rho}
\rho_{0k}=\sum_{j=1}^M\alpha_j\delta_{\vect y_{n_j}\vect y_k},\,\, k=1,\ldots,K,
\end{equation}
where $\delta_{\cdot\cdot}$ is the classical Kronecker delta. Using \eqref{GreenFuncVec} and \eqref{eq:rho}, 
we can write the response matrix  as sum of outer products as follows,
\begin{equation}
\label{responsematrix}
\vect\wP(\omega) \equiv [\wP(\vect x_r,\vect x_s,\omega)]_{r,s=1}^{N}=\sum_{j=1}^M\alpha_j\vect\wg_0(\vect y_{n_j},\omega)\vect\wg_0^T(\vect y_{n_j},\omega)
=\sum_{j=1}^K\rho_{0j}\vect\wg_0(\vect y_{n_j},\omega)\vect\wg_0^T(\vect y_{n_j},\omega).
\end{equation}
Using \eqref{GreenFuncVec}, we also define the $N\times K$ sensing matrix $\vect{\Gc}_0$ as
\begin{equation}\label{sensingmatrix}
\vect{\Gc}_0=[\vect\wg_0(\vect y_1)\,\cdots\,\vect\wg_0(\vect y_K)] \,,
\end{equation}
and write \eqref{responsematrix} in matrix form as
\begin{equation}\label{responsematrix1}
\vect\wP(\omega)=\vect\Gc_0\hbox{diag}({\bf rho_0})\vect\Gc_0^T.
\end{equation}
We  note that the full response matrix $\vect\wP(\omega)$ 
is  symmetric due to Lorentz reciprocity.

Given an array imaging configuration, all the information for imaging is contained in the full response matrix $\vect\wP(\omega)$, including phases. In this case, given a set of illuminations $\{\vect\wf^{(j)}(\omega)\}_{j=1,2,\dots}$, the imaging problem is to determine the location and reflectivities of the scatterers from  the data 
\begin{equation}\label{data}
\vect b^{(j)}(\omega)=\vect\wP(\omega)\vect\wf^{(j)}(\omega)\,\, , \quad j=1,2,\dots 
\end{equation}
received on the array. The components of illumination vectors $\vect\wf^{(j)}(\omega)=[\wf_{1}^{(j)}(\omega),\ldots,\wf_{N}^{(j)}(\omega)]^T$ in \eqref{data} are the signals $\wf_{1}^{(j)}(\omega),\ldots,\wf_{N}^{(j)}(\omega)$ sent from each of the $N$ transducers in the array. 

If only the intensities of the signals are available, the imaging problem is to determine the location and reflectivities of the scatterers from  the absolute value of each component in \eqref{data}, i.e., from the intensity vectors
\begin{equation}\label{dataI}
\vect b_I^{(j)}(\omega)=\hbox{diag}((\vect\wP(\omega)\vect\wf^{(j)}(\omega))(\vect\wP(\omega)\vect\wf^{(j)}(\omega))^*)\,\, , \quad j=1,2,\dots \, .
\end{equation}
In \eqref{dataI}, the superscript $*$ denotes conjugate transpose.
This problem is, however, nonlinear and, therefore, there is much interest in finding algorithms that give the true global solution effectively. 

\section{The time reversal operator} 
\label{sec:timereversalop}

In this paper, we propose a novel imaging strategy for the case in which only data of the form \eqref{dataI}
is recorded and known. The main idea behind the approach proposed here is that we can use a related
matrix to the full response matrix 
$\vect\wP(\omega)$ that has good properties for imaging and can be obtained from data of the form \eqref{dataI}. This related 
matrix is the {\em time reversal matrix} $\vect\wM(\omega)=\vect\wP^*(\omega) \vect\wP(\omega)$. In this Section, 
we will show first how to obtain it from the intensity vectors \eqref{dataI} using the polarization identity, 
and how to use it for imaging using its singular value decomposition.


\subsection{Evaluation of the time reversal operator from quadratic measurements} 
\label{sec:eval_tro}

The key point in active array imaging is that we control the illuminations that probe the medium and, therefore,
we can design illumination strategies favorable for imaging. In our case, we seek an illumination strategy from which
 can obtain the {\em time reversal matrix}  
$\vect\wM(\omega)=\vect\wP^*(\omega) \vect\wP(\omega)$ from \eqref{dataI}. 
Suppose we can put any illumination $\vect\wf(\omega)$ on the array, but we can only measure quadratic measurements as in \eqref{dataI}, i.e., only the intensity of the data can be recorded. 
In that case, we also have access to the quadratic form
\begin{equation}\label{M}
\langle \vect\wf(\omega), \vect\wM(\omega) \vect\wf(\omega) \rangle,~\vect\wM(\omega)= \vect\wP^*(\omega)\vect\wP(\omega)\, .
\end{equation}
Indeed,
\begin{equation}\label{Mbis}
\langle \vect\wf(\omega), \vect\wM(\omega) \vect\wf(\omega) \rangle= \langle \vect\wf(\omega), \vect\wP^*(\omega)\vect\wP(\omega) \vect\wf(\omega) \rangle 
= \langle  \vect\wP(\omega) \vect\wf(\omega), \vect\wP(\omega) \vect\wf(\omega)\rangle =\norm{\vect\wP(\omega) \vect\wf(\omega)}^2.
\end{equation}
Note that only the {\em total power} 
\begin{equation}\label{power}
\norm{\vect\wP(\omega) \vect\wf(\omega)}^2 = \sum_{i=1}^N \abs{\vect\wP(\omega) \vect\wf(\omega)}_i^2
\end{equation}
received at the array is involved in \eqref{Mbis}. In \eqref{power},  $\abs{\vect\wP(\omega) \vect\wf(\omega)}_i^2$ is the intensity
of the signal received at the i-th transducer. Note that $\vect\wM(\omega)$ represents a self-adjoint transformation from the {\em illumination space} $\mC^N$ to the {\em illumination space} $\mC^N$. The entries of this
$N \times N$ square matrix  can be obtained  from the total power received at the array using multiple illuminations as follows.  

The i-th entry in the diagonal $M_{ii}(\omega)$, $i=1,\dots,N$, is just the total power received at the array when only the i-th transducer of the array fires a signal. In other words,
$M_{ii}=\norm{\vect\wP(\omega) \vect\we_i(\omega)}^2$, where
the illumination vector $\vect\we_i = [0,0,\ldots,1,0,\ldots,0]^T$ is the vector whose entries are all
zero except the i-th entry which is $1$.


The off-diagonal terms $M_{ij}(\omega)$, $i \neq j$ can be found from the polarization identity in the complex-valued inner product spaces. Namely, using the polarization identity
\begin{equation}\label{polarization}
2 \langle \vect x, \vect y \rangle = \norm{\vect x+\vect y}^2 - \norm{\vect x}^2 -\norm{\vect y}^2
+ {\bf i} \left(  \norm{\vect x-i \vect y}^2-  \norm{\vect x}^2 -  \norm{\vect y}^2 \right),
\end{equation}
we obtain
\begin{equation}\label{polarization_re}
\mbox{Re}(M_{ij}(\omega))= \mbox{Re}(M_{ji}(\omega)) =  \frac{1}{2} \left( \norm{\vect\wP(\omega) \vect\we_{i+j}}^2  - \norm{\vect\wP(\omega) \vect\we_i}^2 - \norm{\vect\wP(\omega) \vect\we_j}^2\right),  
\end{equation}
using the illumination vector $\vect\we_{i+j}= \vect\we_{i}+\vect\we_{j}$, and
\begin{equation}\label{polarization_im}
\mbox{Im}(M_{ij}(\omega))= - \mbox{Im}(M_{ji}(\omega)) =  \frac{1}{2} \left( \norm{\vect\wP(\omega) \vect\we_{i-{\bf i} j}}^2  - \norm{\vect\wP(\omega) \vect\we_i}^2 - \norm{\vect\wP(\omega) \vect\we_j}^2\right),  
\end{equation}
using the illumination vector $\vect\we_{i-{\bf i} j}= \vect\we_{i}-{\bf i} \vect\we_{j}$. In \eqref{polarization_re} and \eqref{polarization_im}, $\mbox{Re}(\cdot)$ and $\mbox{Im}(\cdot)$ denote the real and imaginary parts
of a complex number, respectively. Again, only the total power received on the array is involved
in these formulas. 

From  \eqref{polarization_re} and \eqref{polarization_im} it follows that we can recover all the entries
in matrix $\vect\wM(\omega)$ by using the following illumination strategy. Send in the illuminations 
$\vect\we_1 =(1,0,0,\dots,0)$, $\vect\we_2 =(0,1,0,\dots,0)$, $\vect\we_{1+2} =(1,1,0,\dots,0)$, and $\vect\we_{1-{\bf i} 2} =(1,- {\bf i},0,\dots,0)$. Then, from the above elementary formulas we can determine the
entries $M_{11}$,  $M_{22}$ and $M_{12}= {\overline M_{21}}$. 
Following the same procedure for each pair of transducers $i$ and $j$ in the array we can determine all four entries $M_{ii}$, $M_{jj}$, $M_{ij}$, and $M_{ji}$.
This means that if measure the total power received at the array from
$N^2$  illuminations we can determine $\vect\wM(\omega)$ completely.

\subsection{Incomplete set of illuminations} 
\label{sec:incompletedata}

In the previous subsection we used the polarization identity to obtain the time reversal matrix  $\vect\wM(\omega)$ using $N^2$  illuminations. In this case, all the entries of the matrix $\vect\wM(\omega)$ can be found. However, there are situations
in which the data from some illuminations are corrupted and must be discarded. In this case, the matrix $\vect\wM(\omega)$  is not full, and the images have to be formed from an incomplete set of illuminations. 
%
We may model these situations by using randomly selected pairs of transducers, and recovering the entries of $\vect\wM(\omega)$ that can not be found from those illuminations by using matrix completion. 
This is possible because  the data matrix $\vect\wM(\omega)$ is of low rank since the image is sparse.
The reconstruction of $\vect\wM(\omega)$ can be accomplished by minimizing its nuclear norm subject to agreement with its known entries. In more detail, we first recover $\vect\wM(\omega)$ by solving the optimization problem
\begin{equation}
\min\|\mathbf{\widehat C} \|_\star \quad\text{s.t.}\quad \widehat C_{ij} = \widehat M_{ij},\, (i,j)\in\Omega,
\label{minnuclear}
\end{equation}
with the singular value thresholding algorithm \cite{Cai08},  and then we apply the two imaging methods proposed in Section~\ref{sec:methods} to the reconstructed matrix 
$\mathbf{\widehat C}$. In \eqref{minnuclear}, $\|\cdot \|_\star$ denotes the nuclear norm of a matrix, and $\Omega$ denotes a random subset 
of $\vect\wM(\omega)$. In \cite{Candes12}, it was proven that most $N \times N$ matrices of rank $r$ can be perfectly recovered from noiseless data by 
solving \eqref{minnuclear}, provided that the cardinality  of $\Omega$ is greater than $c \, N^{6/5} r \log N$, for some constant $c$. Our numerical experiments in Section \ref{sec:numerics} show that we can recover the
noiseless signal when no more than $50\%$ of entries of $\vect\wM(\omega)$ are missing.

Another interesting  intensity-only imaging situation with an incomplete set of illuminations is when one has access to reliable data 
but wants to minimize their number. In this case, one can form the images from data obtained from a few good illuminations. 
The key point here is that the illumination done from the sources at the edges of the array are optimal in the sense that
they carry most of the information needed for imaging~\cite{BPV08}.  Note that, in this situation, the entries of  the data matrix 
$\vect\wM(\omega)$ are not selected at random, 
and matrix completion cannot be accomplished because many rows and columns of $\vect\wM(\omega)$ are unsampled. When the illumination
is done using only a few sources at the edges of the array, only the submatrices at the four corners of $\vect\wM(\omega)$ are known.
Our numerical experiments in Section \ref{sec:numerics} show that intensity-only imaging can be carried out with this partial knowledge of $\vect\wM(\omega)$ directly, that is, without matrix completion. Furthermore, the numerical experiments show that if the data quality is good, i.e., if the signal to noise ratio is high, the number of illuminations needed for imaging can be quite small.

\subsection{The singular value decomposition of the time reversal operator}
\label{sec:sdv_tro}

In this section, we describe how to use two well known imaging methods to obtain images from intensity measurements.
We use an optimization-based method and a subspace projection method. In both methods, we exploit the fact that the SVD of the {\em time reversal matrix}  $\vect\wM(\omega)=\vect\wP^*(\omega) \vect\wP(\omega)$
is similar to the SVD of the full data matrix $\vect\wP(\omega)$, which also contains the information about the phases
of the signals received at the array. Indeed, if we write the SVD of $\vect\wP(\omega)$ in the form
\begin{equation}
\label{svd1}
\vect\wP(\omega)=\vect\wU(\omega)\vect\Sigma(\omega)\vect\wV^\ast(\omega)
=\sum_{j=1}^{\tilde M}\sigma_j(\omega)\wU_j(\omega)\wV_j^\ast(\omega)\, ,
\end{equation}
it follows from the definition of $\vect\wM(\omega)$ \eqref{M}  that the SVD of $\vect\wM(\omega)$ is given by
\begin{equation}
\vect\wM(\omega)= \vect\wV(\omega)\vect\Sigma^2(\omega)\vect\wV^\ast(\omega)\, 
=\sum_{j=1}^{\tilde M}\sigma^2_j(\omega)\wV_j(\omega)\wV_j^\ast(\omega)\, .
\end{equation}
In these equations, $\sigma_1(\omega)\ge\cdots\ge\sigma_{\tilde M}(\omega)>0$ are the nonzero singular values, 
and $\wU_j(\omega)$, $\wV_j(\omega)$ are the corresponding left and right singular vectors, respectively. They fulfill the following equations:
\begin{equation}
\label{svd2}
\vect\wP^\ast(\omega) \wU_j(\omega) = \sigma_j(\omega) \wV_j(\omega)\, , \quad 
\vect\wP(\omega) \wV_j(\omega) = \sigma_j(\omega) \wU_j(\omega)\,,\,\,j=1,\ldots,N .
\end{equation}
Since $\vect\wP(\omega)$ is complex-valued but symmetric, $\wU_j(\omega)= e^{{\bf i} \theta_j} {\overline \wV_j(\omega)}$  for some  unknown global phase $\theta_j$, $j=1,\ldots,N$. Hence, it follows from \eqref{svd2} that
\begin{equation}
\label{svd3}
\vect\wP(\omega) \wV_j(\omega) = \sigma_j(\omega) e^{{\bf i \theta_j}} \overline{\wV_j(\omega)}\,,\,\,j=1,\ldots,N ,
\end{equation}
for an unknown global phase $e^{{\bf i \theta_j}}$ which is different for each singular vector $\wV_j(\omega)$.
Formula~\eqref{svd3} implies that if the singular vector $\wV_j(\omega)$ is the illumination used at the array, 
then the data on the array is known 
up to a global phase. This observation is the key point for the proposed optimization-based algorithm described
in Section \ref{sec:methods}.

Subspace projection algorithms requires another observation. Namely, the matrices $\vect\wP(\omega)$ and
$\vect\wM(\omega)$ have the same kernel. Then, it immediately implies that subspace projection algorithms, e.g.
MUSIC type algorithms, can be applied  to find the locations of the scatterers if the matrix $\vect\wM(\omega)$ is known.

Note that if $\vect\wM(\omega)=\vect\wV(\omega) \vect\Sigma^2(\omega)\vect\wV(\omega)^*$ has been obtained, then 
$\vect\wP(\omega)$  
is the complex-valued symmetric matrix of the form $\vect\wP(\omega)=\overline{\vect\wV(\omega)} D\vect \Sigma(\omega) \vect\wV(\omega)^*$, where  $D$ is an unknown  diagonal matrix with  $e^{{\bf i \theta_k}}$ on the kth diagonal entry. Thus, 
the problem of imaging from intensity-only measurements can be reduced to one in which the full data at the array is known, 
as it is explained next. 
 
\subsection{Sensitivity to noise} 
\label{sec:noise}
Robustness to noise of the proposed approach is a consequence of the central limit theorem, and the fact that we measure 
the total power~\eqref{power} to construct the $N\times N$ {\em time-reversal matrix} $\vect\wM(\omega)$. Indeed, suppose the 
noise at the i-th receiver is modeled by adding a  random variable $\zeta_i$
uniformly distributed  on $[(1-\eps) b_{Ii} , (1+\eps) b_{Ii}]$, 
where $b_{Ii}=\abs{\vect\wP(\omega) \vect\wf(\omega)}_i^2$ is the noiseless intensity
received on the i-th receiver, and $\eps\in (0,1)$ is a parameter that measures the noise strength. If we define the signal-to-noise ratio  at the i-th receiver ($\mbox{SNR}_{i}$)
as the mean to standard deviation of the received power,
then the $\mbox{SNR}_{i}$ on each receiver is the same, and is given by
\[
SNR_i= \frac{b_{Ii}}{\sqrt{Var\left(\zeta_i \right)}}= \frac{\sqrt{3}}{\eps}.
\]
Therefore, the signal-to-noise ratio for the total power is
\[
SNR =\frac{\sum_{i=1}^N b_{Ii}}{\sqrt{\sum_{i=1}^N Var\left(\zeta_i \right)}}= \frac{\sqrt{3}}{\eps} \frac{\sum_{i=1}^N b_{Ii}}{ \sqrt{\sum_{i=1}^N b_{Ii}^2}} \sim 
 O(\sqrt{N}/\eps),
\]
if the intensity does not vary too dramatically from one receiver to another. For example, it suffices to assume there exists $C>0$ so that
\[
\max_i b_{Ii} \leq C \min_i b_{Ii}.
\]

It is straightforward to see that if the intensity at the i-th receiver is a random variable uniformly distributed on
$[(1-\eps) b_{Ii} , (1+\eps) b_{Ii}]$, then the noise in each entry of the time reversal matrix 
$\vect\wM(\omega)$ is a family of zero-mean, uncorrelated Gaussian random variables with variance
\[
\sigma^2=\delta\|\vect\wM(\omega)\|_F^2/N^2\,.
\]
Here, $\|\cdot \|_F$ is the Frobenius matrix norm, and the positive constant is given by
\[
\delta = O\left(\eps^2/N \right).
\] 
Hence, the larger the number of transducers $N$ in the array, the smaller the noise in the resulting {\em time reversal
matrix} used for imaging.

\section{Methods for array imaging}
\label{sec:methods}

In this section, we describe the two imaging methods we use to form the images. At the beginning of each subsection we will assume that the full data matrix $\vect\wP(\omega)$ is recorded and known, i.e., that the amplitudes and the phases of the signals received at the array are available for imaging. At the end of each subsection we show how these methods can be applied to the {\em time reversal matrix} $\vect\wM(\omega)$. 

\subsection{Multiple Measurement Vector imaging method}
We now describe an optimization-based imaging method that exploits the sparsity of the scatterers in the IW.
We will formulate active array imaging as a joint sparsity recovery problem  where we seek an unknown matrix whose 
columns share the same support but possibly different nonzero values. This is known as the
Multiple Measurement Vector (MMV) approach  that has been widely studied in passive source localization  
\cite{MCW05} and active array imaging problems with non negligible multiple scattering \cite{Chai14} with success.
This method can recover the location and reflectivity of the scatterers  exactly from full data in the noise-free case,
and is robust with respect to noise (see \cite{Chai14} for details). Next, we briefly describe the MMV approach assuming that the full data matrix $\vect\wP(\omega)$ is known.

Assume  that the number of scatterers $M$ is much smaller than the number of grid points $K$,  so $M\ll K$.
Hence, the {\em reflectivity vector} $\bfrho_0 = (\rho_{01},\rho_{02},\dots,\rho_{0K})\in\mC^K$, is sparse.
From \eqref{data}, the signal scattered back from the scatterers and received on the array is given by 
$\vect\wP(\omega)\vect\wf(\omega)$, where $\vect\wf(\omega)$ is the illumination sent from the array. 
Then,  we can define the
linear operator $A_{\wf(\omega)}$ that relates the {\em reflectivity vector} $\bfrho_0$ with the received signals through the identity
\begin{eqnarray}
\vect\wP(\omega)\vect\wf(\omega)
&=&\sum_{j=1}^K\rho_{0j}(\vect\wg_0^T(\vect y_j,\omega)\vect\wf(\omega))\vect\wg_0(\vect y_j,\omega)
=A_{\widehat f(\omega)}\bfrho_0. \label{mop1}
\end{eqnarray}
Hence, $A_{\widehat f(\omega)}$ is  the $N\times K$ matrix 
\begin{equation}
\label{eq:A_matrix}
A_{\wf(\omega)}= 
\begin{bmatrix}
\hat g_{\hat f}(\vect y_1,\omega)\vect\wg(\vect y_1,\omega) & \hat g_{\hat f}(\vect y_2,\omega)\vect\wg(\vect y_2,\omega)& \cdots & \hat g_{\hat f}(\vect y_K,\omega)\vect\wg(\vect y_K,\omega) \,  \\
\end{bmatrix} \, 
\end{equation}
that depends on the illumination. In \eqref{eq:A_matrix}, $\wg_{\wf}(\vect y_j,\omega)=\vect\wg^T_{0}(\vect y_j,\omega)\vect\wf(\omega)$, $j=1,\dots,K$, are scalars that represent the field at $\vect y_j$ due to the illumination $\vect\wf(\omega)$ sent from the array. With this notation,
active array imaging with a single illumination amounts to solving for $\bfrho_0$ from the system of equations 
\begin{equation}\label{eq:single}
\vect\Ac_{\wf(\omega)} \bfrho_0=\vect b(\omega).
\end{equation}
Since the number of transducers $N \ll K$ in the IW, the system of equations \eqref{eq:single} is underdetermined and, therefore, 
there are many configurations of scatterers that match the data vector $\vect b(\omega)$. 
However, due to the known sparsity of the reflectivity vector $\bfrho_0$, one can use $\ell_1$ minimization 
\begin{equation}\label{l1singleillum}
\min\|\bfrho\|_{\ell_1}\quad\quad\text{s.t.}\quad\vect\Ac_{\wf(\omega)} \bfrho=\vect b(\omega)\, ,
\end{equation}
to find the sparsest solution from noiseless data. It is well known that under certain conditions on 
the operator $\Ac_{\wf(\omega)}$, and on the sparsity of $\bfrho_0$,  $\ell_1$ minimization is equivalent to $\ell_0$ minimization \cite{Candes05,Donoho06}.
When the data $\vect b(\omega)$ is contaminated by a noise vector $\vect e$,
then one can solve the relaxed problem
\begin{equation}\label{l1singleillumnoise}
\min\|\bfrho\|_{\ell_1}\quad\quad\text{s.t.}\quad\|\vect\Ac_{\wf(\omega)} \bfrho-\vect b(\omega)\|_{\ell_2}<\varepsilon \, ,
\end{equation}
for some given positive constant $\varepsilon$. 
The full data vector $\vect b(\omega)$ in \eqref{eq:single}-\eqref{l1singleillumnoise}, which contains both
the amplitudes and the phases of the collected signals, is obtained from a
single illumination $\vect\wf(\omega)$.

When multiple illuminations are available, one could solve the $\ell_1$ minimization problem
\begin{equation}\label{l1normminimization}
\min\|\bfrho\|_{\ell_1}\quad\text{s.t. }\quad
\|\Ac_{\wf^{(j)}(\omega)}\bfrho-\vect b^{(j)}(\omega)\|_{\ell_2}\le\varepsilon\quad
\mbox{for}\quad
j=1,2,\ldots,\nu 
\end{equation}
to capture the sparsity of $\bfrho_0$. Here, $\nu$ is the number of illuminations. This formulation, however,
does not exploit the data structure optimally, as the solution vectors from different illuminations have the same
support. To take advantage of the data structure, one can formulate the problem of array imaging with multiple illuminations as a joint sparse recovery problem, also known as the MMV formulation, aims to recover unknown sparse matrices with nonzero entries restricted to a small number of rows \cite{Cotter05,MCW05,CH06,ER10}.

We use this formulation for active array imaging in two steps as in \cite{Chai14}. In the first step, we determine the locations of the scatterers that are treated as equivalent sources. The equivalent sources have unknown locations but strengths related, in a known way, to the reflectivities of the scatterers and to the used illuminations. In the second step, once the locations of the scatterers have been obtained, we recover the true reflectivities  easily from these known relationships.

\subsubsection{Locations of the scatteters}

In the first step, the sought matrix  is the ${K \times \nu}$ matrix 
$\vect{\Rho}_0=[{\bfgamma}_0^{(1)}\,\ldots\,{\bfgamma}_0^{(\nu)}]$ whose
$j^\mathrm{th}$ column corresponds to the {\em effective source vector} ${\bfgamma}_0^{(j)}$ whose components are given by
\begin{equation}\label{effsource}
(\bfgamma_0^{(j)})_k = \wg_{\wf^{(j)}}(\vect y_k,\omega) \rho_{0k}\, , \,\, k=1,\dots,K,
\end{equation}
under illumination $\vect\wf^{(j)}(\omega)$, $j=1,\ldots,\nu$. This matrix variable $\vect\Rho_0\in\mC^{K \times \nu}$ has columns that share the same sparse support but possibly have different nonzero values due to the different illuminations.

The MMV formulation for active array imaging is to solve for $\vect\Rho_0$ from the matrix-matrix equation
\begin{equation}\label{linearsystemmmv}
\vect\Gc_0 \vect\Rho = \mathbf{B},
\end{equation}
where $\vect{\Gc}_0$ is the $N\times K$ sensing matrix \eqref{sensingmatrix}, and  $\mathbf{B}=[\vect b^{(1)}\,\ldots\,\vect b^{(\nu)}]$ is 
the ${N \times \nu}$ data matrix whose columns are the full data vectors generated by the $\nu$ illuminations. 
In the MMV framework, the sparsity of the matrix variable $\vect\Rho$ is characterized by the number of nonzero
rows, i.e., by the row-wise $\ell_0$ norm of $\vect\Rho$. 
More precisely, we define the row-support of
a given matrix $\vect\Rho$ by
\[\operatorname{rowsupp}(\vect{\Rho})=\{i:\,\,\|X_{i\cdot}\|_{\ell_2}\neq0\},\]
so the  sparsity of $\vect\Rho$ is measured as $\Xi_0(\vect\Rho)=|\operatorname{rowsupp}(\vect\Rho)|$.
Here, the $i^{th}$ row of $\vect\Rho$ is denoted by $X_{i\cdot}$.
With these definitions, the sparsest solution to \eqref{linearsystemmmv} is given by
\begin{equation}\label{MMV.NP}
\min\Xi_0(\vect\Rho)\quad\text{s.t.}\quad\vect\Gc_0\vect\Rho =\mathbf{B}.
\end{equation}
Since \eqref{MMV.NP}  is an NP hard problem, we solve instead the convex relaxed problem
\begin{equation}\label{MMV21}
\min J_{2,1}(\vect\Rho)\quad\text{s.t.}\quad\vect\Gc_0\vect\Rho=\mathbf{B},
\end{equation}
with the $(p,q)$-norm function $J_{p,q}(\cdot)$ defined as
\begin{equation}
J_{p,q}(\mathbf{Y})=\left(\sum_{i=1}^m\|Y_{i\cdot}\|_{\ell_p}^q\right)^{1/q}.
\label{eq:Jpq}
\end{equation}
The $(p,q)$-norm function \eqref{eq:Jpq} is simply the $\ell_q$ norm of the vector formed by the $\ell_p$ norms of all the rows of a matrix.

When the  data is contaminated by additive noise vectors $\vect{e}^{(j)}$, $j=1,\ldots,\nu$, 
we solve
\begin{equation}\label{MMV21noise}
\min J_{2,1}(\vect\Rho)\quad\text{s.t.}\quad\|\vect\Gc_0\vect\Rho-\mathbf{B}\|_F<\varepsilon \, ,
\end{equation}
for some pre-specified constant $\varepsilon$.

Formulations \eqref{MMV21} and  \eqref{MMV21noise} have been studied thoroughly during the last few years, see for example \cite{Cotter05,MCW05,CH06,ER10, Chai14}.  Under certain conditions on the matrix $\vect\Gc_0$ and the sparsity of $\vect\Rho_0$, \eqref{MMV21}  recovers the sparsest solution exactly if the data is noise-free. If
 the data is contaminated by additive noise, then \eqref{MMV21noise} recovers the sparsest solution upon a certain error bound. See \cite{Chai14} for more details.

\subsubsection{Using MMV with intensity-only measurements}
It follows from the discussion in Section \ref{sec:timereversalop} that the active array imaging problem with intensity-only measurements can be solved from the knowledge of $\vect\wM(\omega)$ using the MMV framework if the data is generated with 
illumination vectors equal to the right singular vectors of $\vect\wM(\omega)$. 
More specifically, we can consider the MMV formulation~\eqref{MMV21} or~\eqref{MMV21noise}
with $\mathbf{B}=[\vect b^1\,\ldots\,\vect b^\nu]$ being the matrix whose columns are the full
data vectors generated by the illuminations $\wV_j(\omega)$  (up to a global phase), that is,
$\vect b^{(j)}=\sigma_j(\omega) \overline{\wV_j(\omega)}$. In ~\eqref{MMV21} and~\eqref{MMV21noise}, $\vect{\Rho}=[{\bfgamma}^1\,\ldots\,{\bfgamma}^\nu]$  is the unknown matrix whose $j^\mathrm{th}$ column 
corresponds to the {\em effective source vector} defined in \eqref{effsource} including a global phase $e^{-{\bf i \theta_j}}$. Then, we can use~\eqref{MMV21} 
or~\eqref{MMV21noise} to find the locations of the effective sources. 

There are different algorithms for solving \eqref{MMV21} and  \eqref{MMV21noise}. We use an extension 
of an iterative algorithm proposed in \cite{MNPR12} for matrix-vector equations. This method, called GeLMA,
is a shrinkage-thresholding algorithm for solving $\ell_1$-minimization problems which has proven to be very efficient and whose solution does not depend on the regularization parameter that promotes sparse solutions,
see \cite{MNPR12} for more details. We summarize it  for MMV problems in Algorithm~\ref{algo} below.
\begin{algorithm}
\begin{algorithmic}
\REQUIRE 
Set $\vect\Rho=\vect0$, $\vect\Zc=\vect 0$, and pick the step size $\beta$ and the regularization parameter $\tau$.
\REPEAT
\STATE Compute the residual $\vect\Rc= \mathbf{B} - \vect\Gc_0\vect\Rho$ 
\STATE $\vect\Rho\Leftarrow\vect\Rho + \beta\vect\Gc_0^\ast(\vect\Zc + \vect\Rc)$
\STATE $\Rho_{i\cdot}\Leftarrow\operatorname{sign}(\|\Rho_{i\cdot}\|_{\ell_2}-\beta\tau)\frac{\|\Rho_{i\cdot}\|_{\ell_2}-\beta\tau}{\|\Rho_{i\cdot}\|_{\ell_2}}\Rho_{i\cdot}$, $i=1,\ldots,K$
\STATE $\vect\Zc\Leftarrow\vect\Zc + \beta\vect\Rc$
\UNTIL{Convergence}
\end{algorithmic}
\caption{GeLMA-MMV for solving \eqref{MMV21} and  \eqref{MMV21noise}.}
\label{algo}
\end{algorithm}

\subsubsection{Reflectivities of the scatteters}
Once we obtain from \eqref{MMV21} or \eqref{MMV21noise} the matrix $\vect\Rho_\star$ whose columns are the effective sources corresponding to the different illuminations, we estimate the  reflectivities easily by using \eqref{effsource}. More precisely, for each component $i$ in the support of the solution given by \eqref{MMV21} or \eqref{MMV21noise},  we compute the estimated reflectivities $\rho^{(j)}_{\star i}$ corresponding to each illumination $j$ 
as 
\begin{equation}\label{rhostar}
\rho^{(j)}_{\star i} = (\bfgamma_\star^{(j)})_i / \wg_{\wf^{(j)}}(\vect y_i,\omega). 
\end{equation}
We then take the average $\frac{1}{\nu}\sum_{j=1}^\nu\rho_{\star i}^{(j)}$ as the estimated reflectivity. We note that if the noise in the data is high, this last step can bring some ghosts to the final image because $\wg_{\wf^{(j)}}(\vect y_i,\omega)$ can be very small at some locations.  Nevertheless, this last step can be easily avoided by a further regularization as, for example, carrying on the division
only at those pixels where $\wg_{\wf^{(j)}}(\vect y_i)$ is above a certain threshold.

\subsection{Multiple signal classification method}
The MUltiple SIgnal Classification method (MUSIC) is a subspace projection algorithm that uses the 
SVD of the full data array response matrix $\vect\wP(\omega)$ to form the images. 
It is a direct algorithm widely used to image the locations of $M<N$ point-like scatterers in a region of interest. 
Once the locations are known, their reflectivities can be found from the recorded intensities using convex optimization as shown below.

\subsubsection{Locations of the scatterers}
The search of the locations  of the $M$ scatterers is the combinatorial part of the imaging problem and, hence, by far the most difficult task. 
Note that $\vect\wP(\omega)$ is a linear transformation from the {\em illumination space} $\mC^N$ to the {\em data space} $\mC^N$. According to \eqref{svd1}, the
 illumination space can be decomposed into the direct sum of
a  signal space, spanned by the principal singular vectors $\wV_j(\omega)$, $j=1\dots,M$, having non-zero singular values,
and a noise space spanned by the singular vectors having zero singular values.
Since the singular vectors
$\widehat{V}_j(\omega)$, $j=M+1,\ldots,N$, span the noise space, the probing vectors $\vect\wg_0(\vect y_j,\omega)$ will be orthogonal to the noise space
only when $\vect y_j$ corresponds to a scatterer location $\vect y_{n_j}$. Hence, it follows that the scatterer locations must correspond to the peaks
of the functional
\begin{equation}
\label{MUSIC_0}
\mathcal{I}(\vect y_s)=\frac{1}{\sum_{j=M+1}^{N} |\vect\wg_0^T(\vect y_s,\omega)\widehat{V}_j(\omega) |^2 },\,\,s=1,\ldots,K.
\end{equation}
We can interpret \eqref{MUSIC_0} in terms of the images created by the singular vectors having zero singular value, as 
$\vect\wg_0^T(\vect y_s,\omega)\widehat{V}_j(\omega)$ is the incident field at the search point $\vect y_s$ due to a illumination vector 
$\widehat{V}_j(\omega)$ on the array. 
According to this interpretation, the singular vectors having zero singular value do not illuminate the scatterers locations and, hence, \eqref{MUSIC_0}
has a peak when $\vect y_s=\vect y_{n_j}$.

Since in our application the number of scatterers is small, the signal space is much smaller than the noise space and, therefore, it is more efficient to compute
the equivalent functional
\begin{equation}
\label{MUSIC}
\mathcal{I}_{MUSIC}(\vect y_s)=\frac{\min_{1\le j\le K}\|\mathcal{P}\vect\wg_0(\vect y_j,\omega)\|_{\ell_2}}{\|\mathcal{P}\vect\wg_0(\vect y_s,\omega)\|_{\ell_2}},\,\,s=1,\ldots,K,
\end{equation}
with the
projection onto the noise space defined as
\begin{equation}
\label{proyection}
\mathcal{P}\vect\wg_0(\vect y,\omega)=\vect\wg_0(\vect y,\omega) -\sum_{j=1}^M (\vect\wg_0^T(\vect y,\omega)\widehat{V}_j(\omega))
\widehat{V}_j(\omega).
\end{equation}
The numerator in \eqref{MUSIC}  is just a normalization. We note that \eqref{MUSIC}  is robust to noise, even for single frequency and for non-homogeneous, random 
media, and it is quite accurate
for large arrays~\cite{BTPB02}. Generalizations of MUSIC for multiple scattering and extended scatterers have also 
been developed (see, for example, \cite{Gruber04} and \cite{hou06}).

\subsubsection{Using MUSIC with intensity-only measurements}

It is an immediate consequence of the discussion in subsection \ref{sec:sdv_tro} that \eqref{MUSIC} can also be used in the case in which the phases of the data are not recorded. Both, $\vect\wM(\omega)$ and $\vect\wP(\omega)$, share the same right singular vectors and, hence, \eqref{MUSIC} can be applied, without any modification, to determine the location of the scatterers,
once the {\em time reversal matrix} $\vect\wM(\omega)$ has been obtained.

\subsubsection{Reflectivities of the scatteters}\label{whatnot}
Once the locations of the scatterers have been found from \eqref{MUSIC}, we may want to estimate their reflectivities in a second step.  This is still a nonlinear problem as
only the intensities are available. To linearize the problem we follow the same approach proposed in \cite{CMP11}, but restricted
to the support of the solution found from \eqref{MUSIC}. Thus, we introduce the positive semidefinite matrix
\begin{equation}
\RhoQ_\star = \rho_\star\rho_\star^* \in\mathbb{R}^{L\times L}\, ,
\end{equation}
associated with the unknown reflectivities $\rho_\star=[\rho_{\star 1},\ldots,\rho_{\star L}]^T\in\mC^L$ defined in the support $\Lambda_\star$ recovered
in the first step.  Note that now  $|\Lambda_\star | = L \ll K$ and, therefore,  $\RhoQ_\star$ has small dimensions. Following \cite{CMP11},
we could obtain $\RhoQ_\star$ from intensity-only measurements by solving
\begin{equation}\label{linearform}
\Lc_{\wf(\omega)}(Y_\star)=\vect b_I(\omega) \, ,
\end{equation}
where $\Lc_{\wf(\omega)}(Y):=\hbox{diag}(\Ac_{\wf(\omega)}Y\Ac_{\wf(\omega)}^\ast)$ is a linear map from
$\mR^{L\times L}$ to $\mR^{N}$.  An estimate for $\RhoQ_\star$ could be found, in principle, by solving \eqref{linearform} 
by least squares. Note, however, that $\RhoQ_\star$ is of low rank (in fact rank $1$ since it is defined via an outer-product), so we obtain $\RhoQ_\star$  from the following affine rank minimization problem
\begin{equation}\label{minimumrank}
\hbox{rank}(X)\quad\text{subject to}\,\, \Lc_{\wf(\omega)}(X)=\vect
b_I(\omega),
\end{equation}
in order to take advantage of the  additional information on the unknown $Y$. Once $Y$ is found from this optimization problem, we can obtain the amplitude of the reflectivities  by  taking $\bfrho=\sqrt{\hbox{diag}(Y)}$ on the support $\Lambda_\star$. 

However, \eqref{minimumrank} is an NP-hard problem and, therefore, there is no simple algorithm that gives the true global solution effectively. Therefore, we replace $\hbox{rank}(X)$ by the nuclear norm $\|X\|_\ast$ in the objective
function of \eqref{minimumrank}, and consider the following optimization problem as given in \cite{CMP11}
\begin{equation}\label{minimumnuclearnorm}
\min\|X\|_\ast\quad\text{subject to  }\,\,
\Lc_{\wf(\omega)}(X)=\vect b_I(\omega).
\end{equation}
The nuclear norm $\|\cdot\|_\ast$ is the sum of the singular values of the matrix while the rank is the number of nonzero singular values and, hence, it can be used as a convex surrogate for the rank functional \cite{Recht10}. Problem \eqref{minimumnuclearnorm} is now convex and can be solved in polynomial time. 

To solve \eqref{minimumnuclearnorm}, we follow~\cite{TY} and use the gradient descent method with singular value thresholding, as outlined below in 
Algorithm~\ref{algominimumlsqrank}. 
In Algorithm~\ref{algominimumlsqrank}, the soft-thresholding operation  is given by
\begin{equation}
S_{\tau}(G) = \vect\wU diag(\vect\sigma - \tau)^+\vect\wV^\ast,
\
\label{eq:softthresh}
\end{equation}
where $\vect\sigma$ is the vector of positive singular values arranged in descending order,  $\tau >0$ is the thresholding parameter, superscript $+$ means 
positive part,  and $\vect\wU$ and $\vect\wV$ are
the orthogonal matrices from the SVD of $G$.
We stress that through step one, i.e. by using MUSIC to locate the scatterers,
we have effectively reduce the dimension of the unknown $X$ in \eqref{minimumnuclearnorm} and, thus, the optimization problem
is very easy to solve.
\begin{algorithm}
\begin{algorithmic}
\REQUIRE 
Set $Y_{-1}=Y_0=0$ and $t_{-1}=t_0=1$, and pick the initial value for step size $\beta$.
\REPEAT
\STATE Compute weight $w=\frac{t_{k-1}-1}{t_k}$.
\STATE Compute $W_k=(1+w)Y_k-w Y_{k-1}$.
\STATE Compute the matrix $G=W-\beta \Lc_{\wf(\omega)}^\ast(\Lc_{\wf(\omega)}(W)-\vect b_I(\omega))$.
\STATE Set $Y_{k+1}=S_{\tau}(G)$.
\STATE Compute $t_{k+1}=\frac{1+\sqrt{1+4t_k^2}}{2}$.
\UNTIL{Convergence}
\end{algorithmic}
\caption{Iterative algorithm for~\eqref{minimumnuclearnorm}}
\label{algominimumlsqrank}
\end{algorithm}

In Algorithm~\ref{algominimumlsqrank}, the adjoint operator $\Lc_{\wf(\omega)}^\ast:\,\mR^{N}\rightarrow\mR_+^{L\times L}$ is given by
\begin{equation}
\Lc_{\wf(\omega)}^\ast(\vect c)=\Ac_{\wf(\omega)}^\ast\operatorname{diag}\big(\vect c\big)\Ac_{\wf(\omega)} \quad 
\mbox{for} \,\,\vect c \in\mR^{N},
 \end{equation}
which is found from the relation $\langle\Lc_{\wf(\omega)}(Y),\vect c\rangle=\langle Y,\Lc_{\wf(\omega)}^\ast(\vect c)\rangle$.

We have seen in our numerical experiments that replacing the soft-thresholding operation by a rank $1$ enforcement, that is, setting 
$Y_{k+1}=\sigma_1\wU_1\wV_1^\ast$ at each iteration in Algorithm~\ref{algominimumlsqrank}, also gives excellent results. This can be understood as solving
the least squares problem with the rank constrain
\begin{equation}\label{minimumlsqrank}
\min\|\Lc_{\wf(\omega)}(Y)-\vect b_I(\omega)\|\quad\text{subject to}\,\, \hbox{rank}(Y)=1.
\end{equation}
This problem is, however, non-convex due to the non-convexity of the set of low-rank matrices and, therefore, it might
not converge to the true solution in general.

\section{Numerical experiments}
\label{sec:numerics}
In this section we present numerical simulations  in two dimensions. The linear array 
consists of $100$ transducers that are one wavelength $\lambda$ apart. 
Scatterers are placed within an IW  of size $30\lambda\times30\lambda$ which is at a distance $L=100\lambda$ from the linear array. 
The amplitudes of the reflectivities of the scatterers and their phases are set randomly in each realization.  
The scatterers are within an IW that is discretized using a uniform lattice with points separated by one wavelength 
$\lambda$. This results in a $30\times30$ uniform mesh. Hence, we have $900$ unknowns.
In all the images shown below, we normalize the spatial units by the wavelength $\lambda$.

\begin{figure}[t]
\begin{center}
\begin{tabular}{ccc}
\includegraphics[scale=0.25]{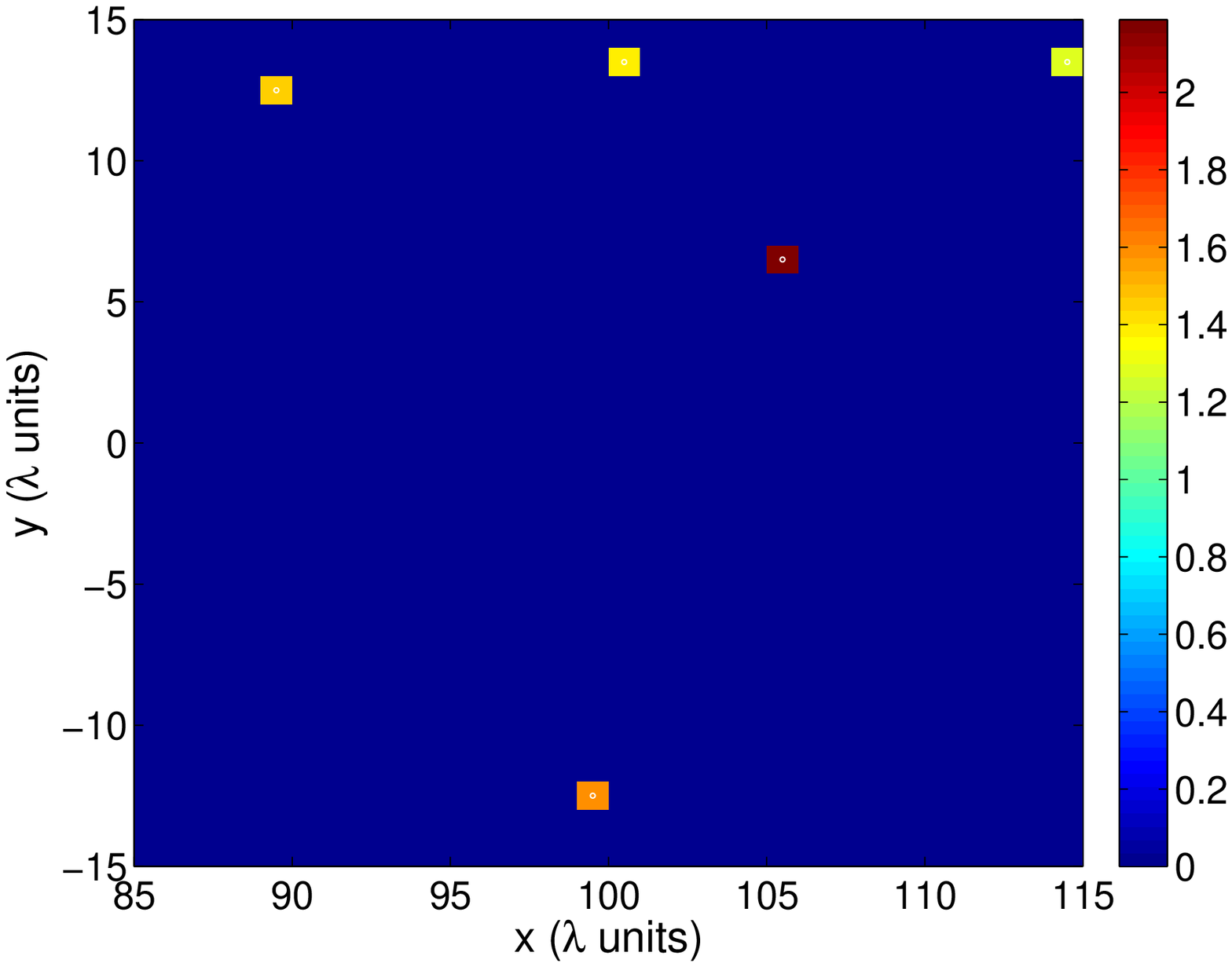} & 
\includegraphics[scale=0.25]{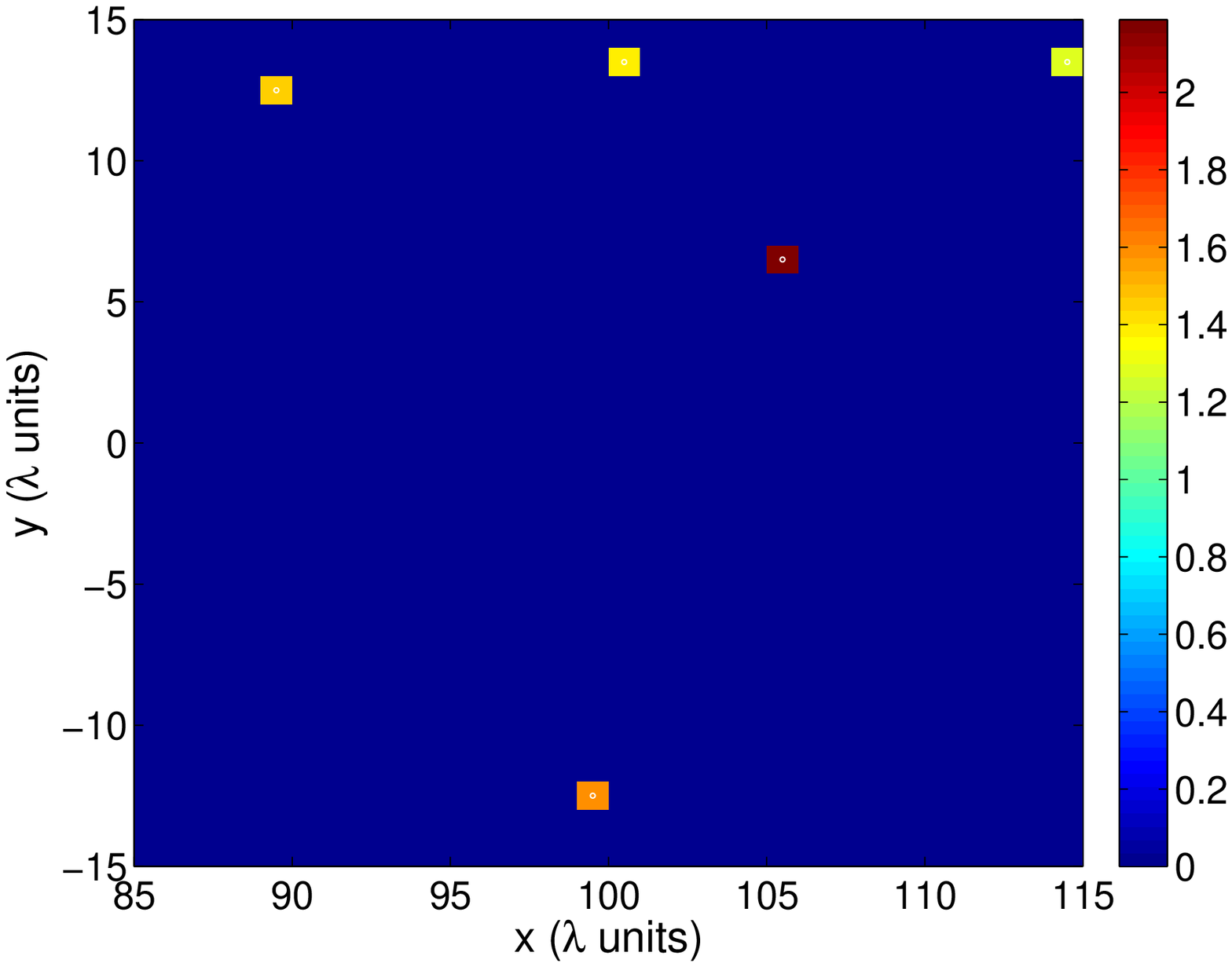} & 
\includegraphics[scale=0.25]{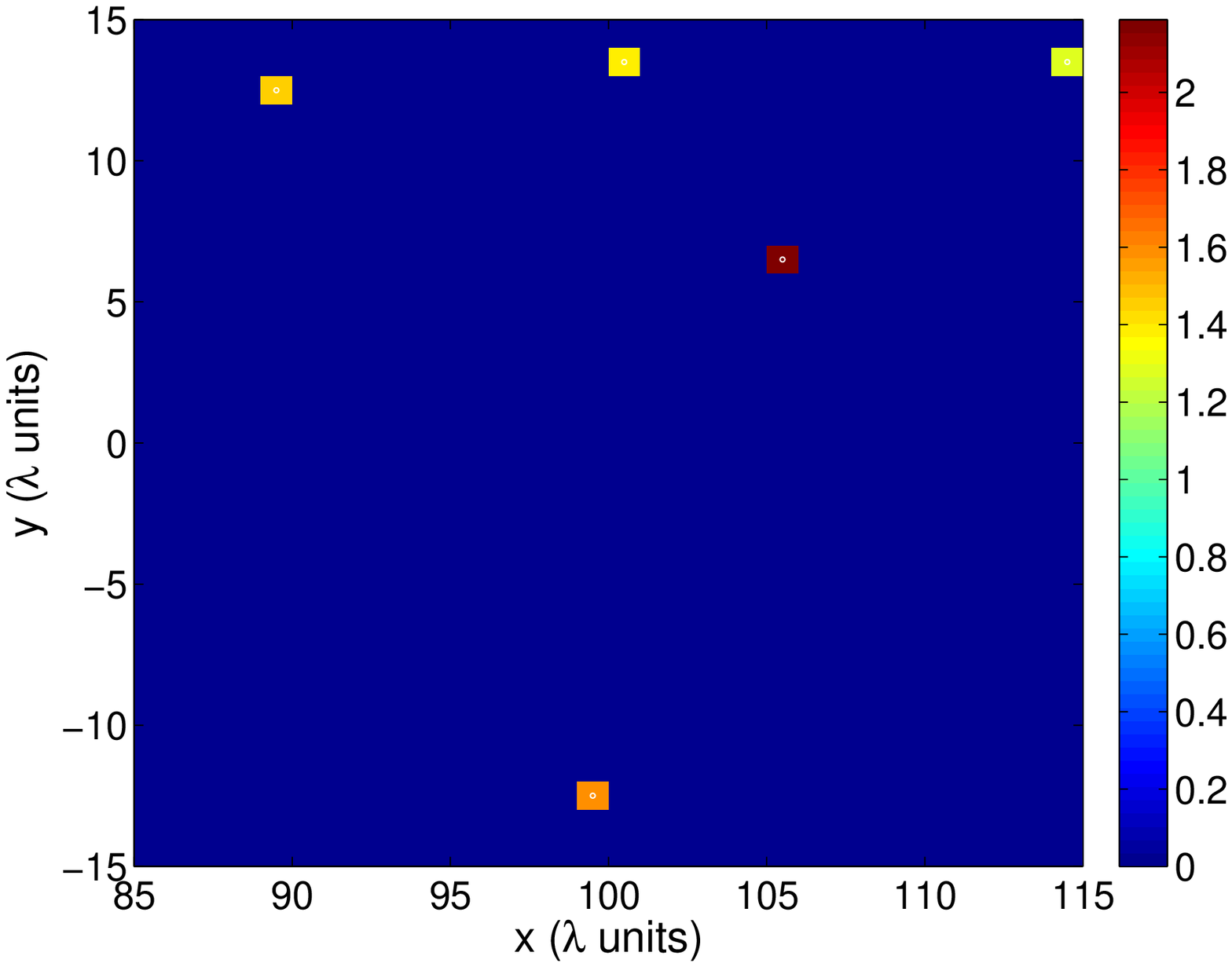} \\
\includegraphics[scale=0.25]{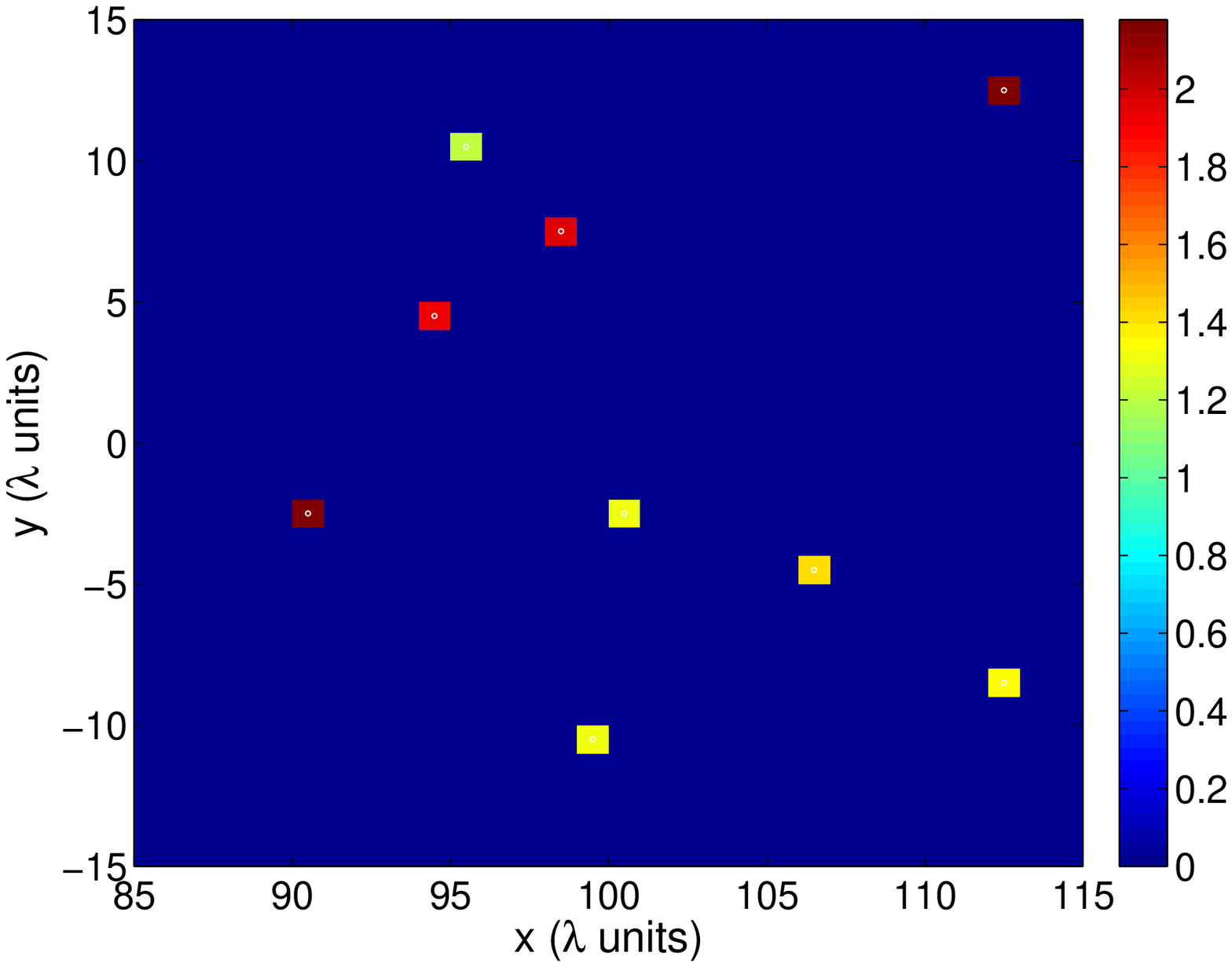} & 
\includegraphics[scale=0.25]{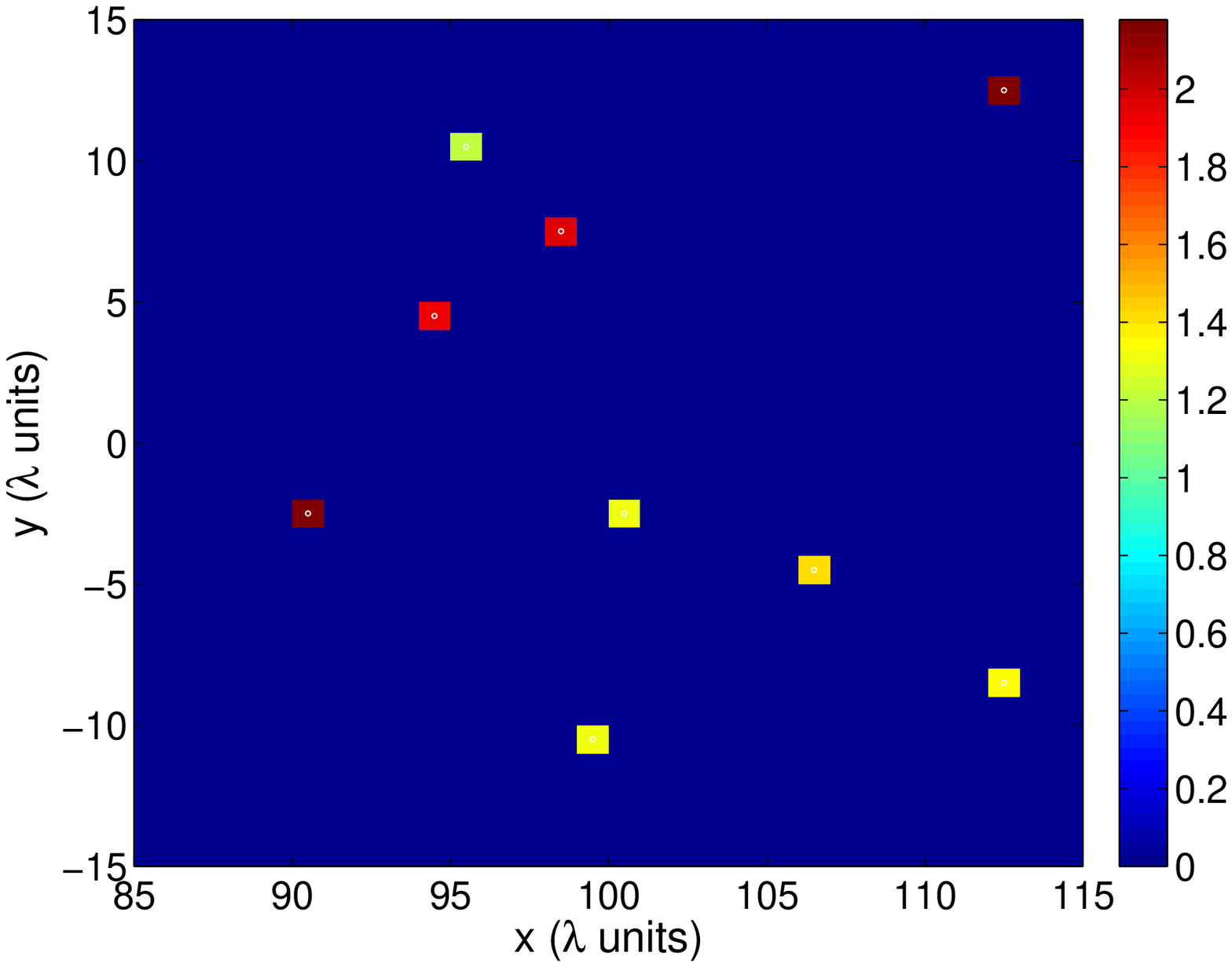} & 
\includegraphics[scale=0.25]{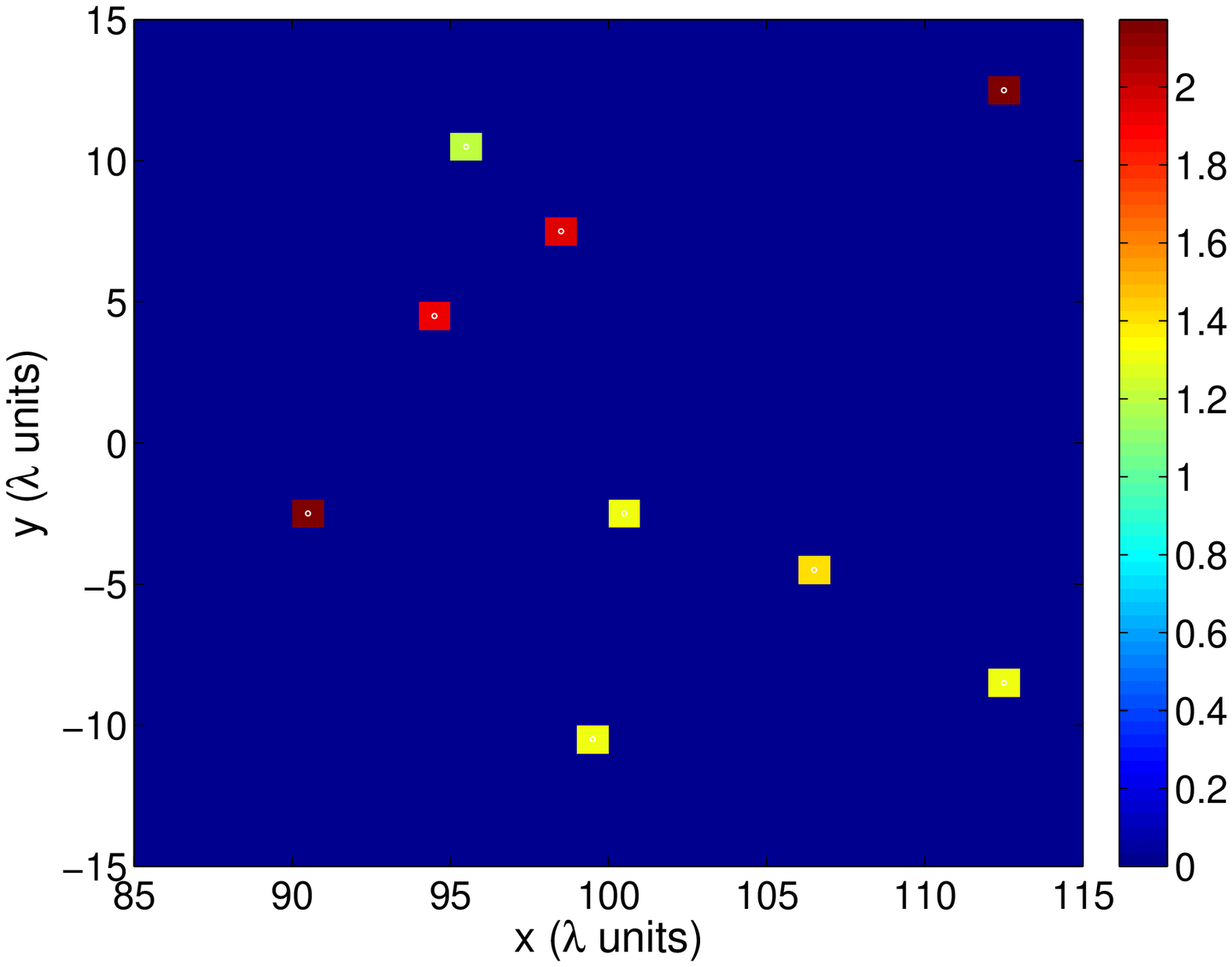} 
\end{tabular}
\caption{Noiseless data. Top and bottom rows are two different configurations  with $5$ and $9$ scatterers, respectively.
The left column shows the  original configurations of the scatterers. The middle and right columns show the amplitudes of the reflectivities 
obtained with MUSIC after nuclear norm minimization, and with the MMV formulation, respectively.
}
\label{fig:nonoise}
\end{center}
\end{figure}

Figure \ref{fig:nonoise} shows the images obtained with MUSIC (middle column) and with the MMV formulation 
(right column) using noisless data. The top and bottom rows are two different configurations with $5$ and $9$ scatterers, respectively. The left column shows the  distribution of scatterers to be recovered.
When there is no noise in the data, both methods recover the positions and reflectivities of the scatterers exactly.
The exact locations of the scatterers in these images are indicated with small white dots.

Next, we examine the performance of these two methods when noise is added to the data.
We simulate instrument noise by adding a random variable uniformly distributed, $\zeta_i$, to the noiseless intensity 
$b_{Ii}^{(j)}=\abs{\vect\wP(\omega) \vect\wf^{(j)}(\omega)}_i^2$ received on each transducer $i$, $i=1,\dots,N$, 
when the vector $\vect\wf^{(j)}$ illuminates the IW (see subsection \ref{sec:noise}). 
With this model, the intensity recorded at the i-th transducer is $[(1-\eps) b_{Ii}^{(j)} , (1+\eps) b_{Ii}^{(j)}]$, where $\eps\in(0,1)$ denotes the strength of the noise.

\begin{figure}[t]
\centering
\begin{tabular}{ccc}
\includegraphics[scale=0.25]{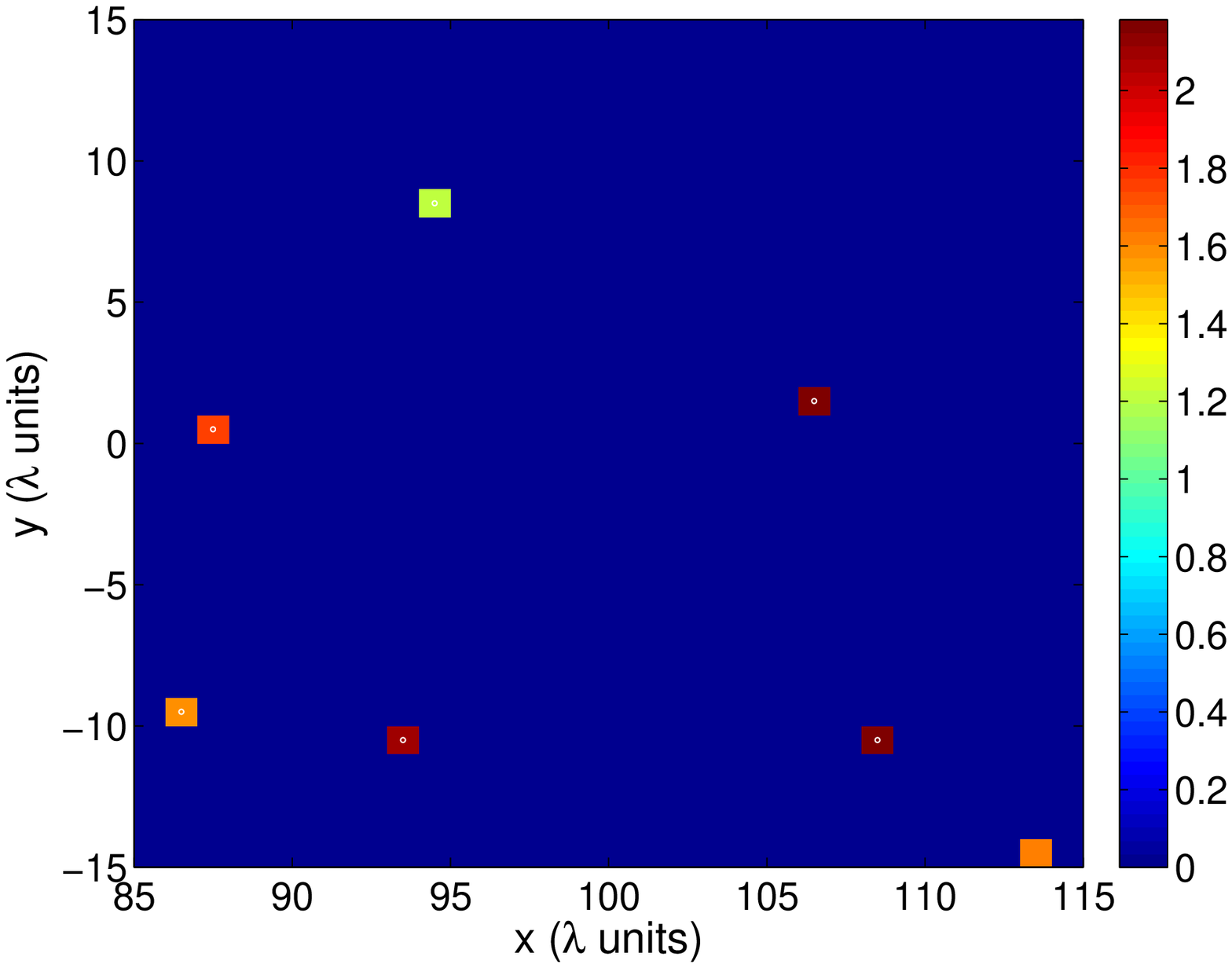} & 
\includegraphics[scale=0.25]{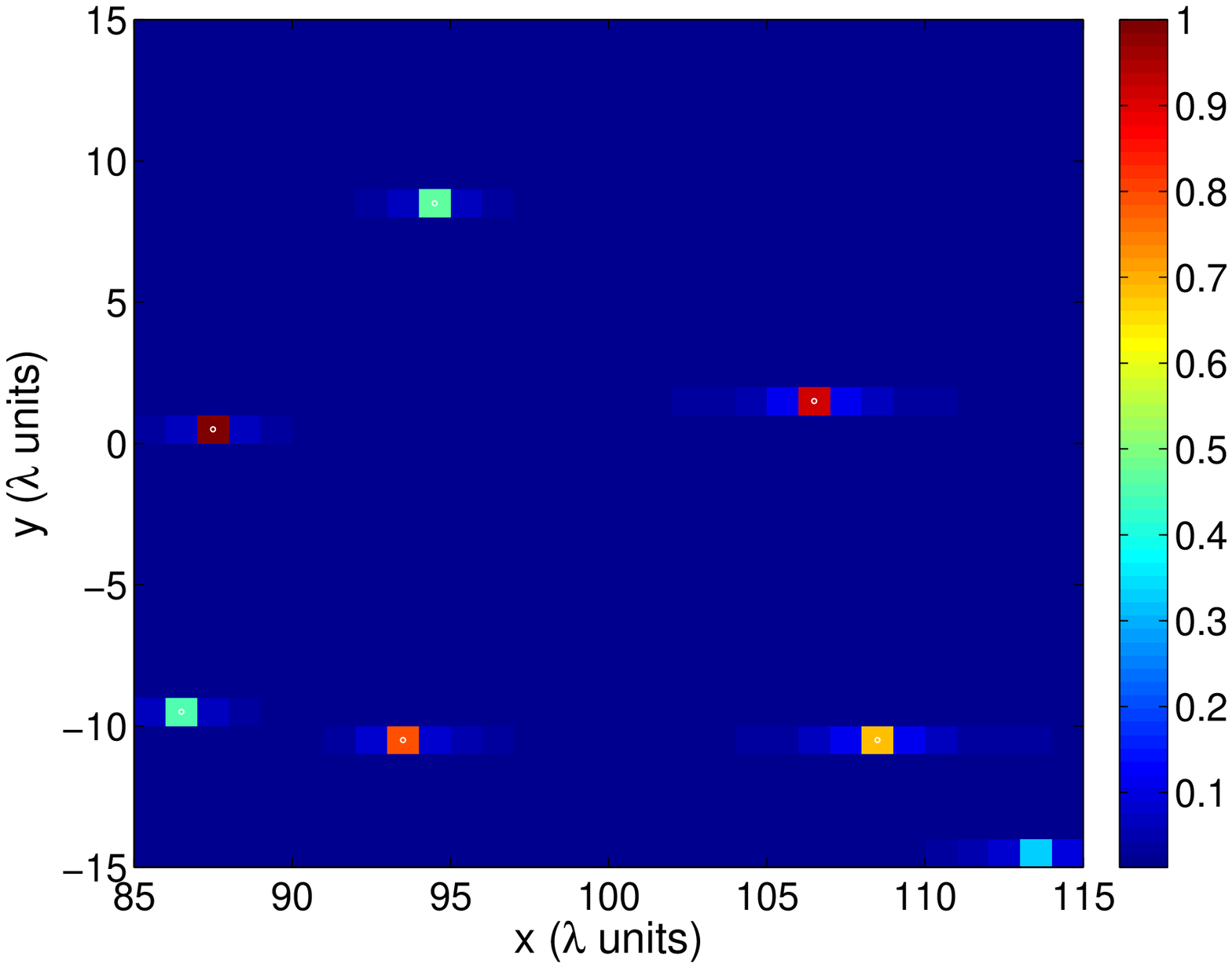} & 
\includegraphics[scale=0.25]{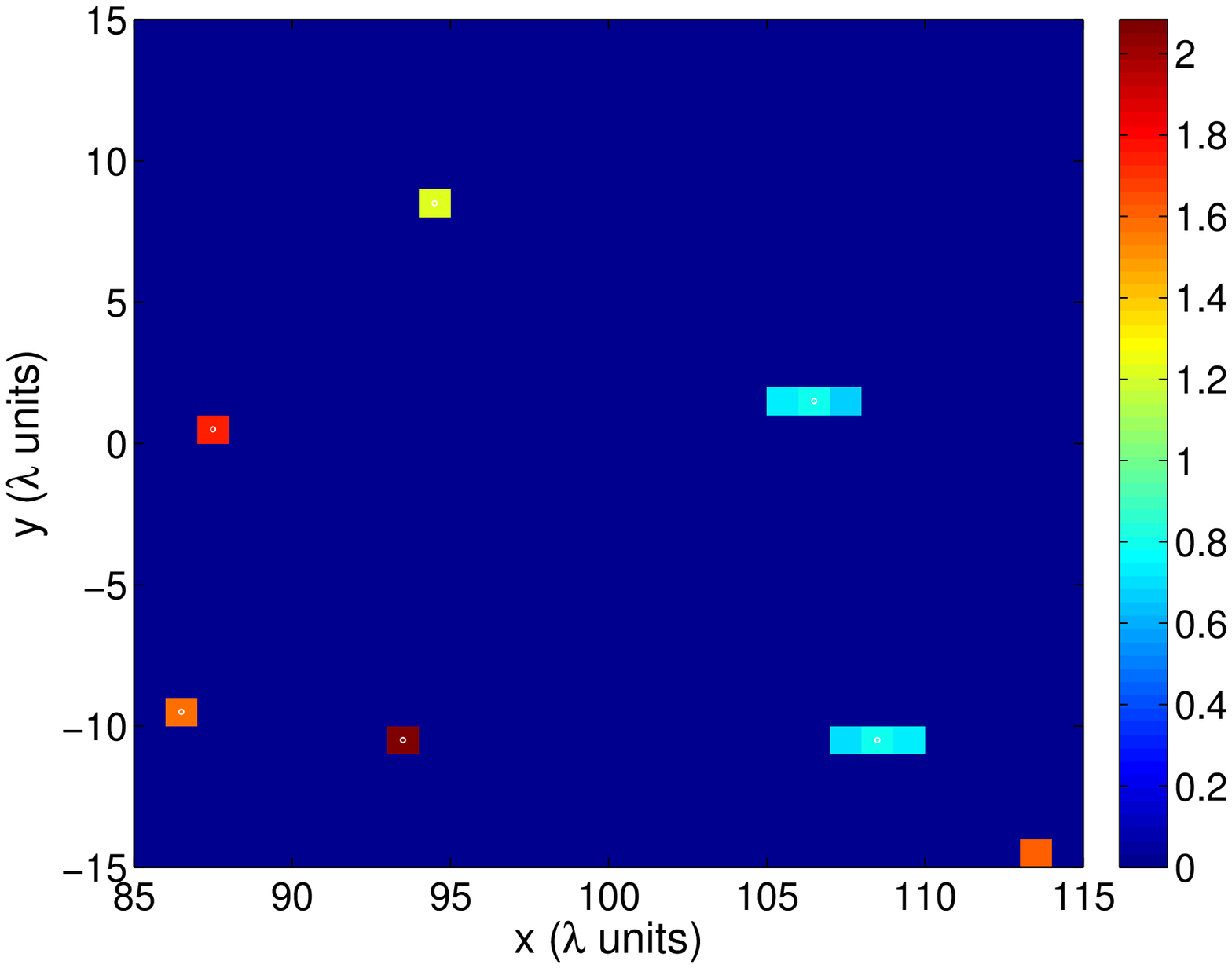} \\
\includegraphics[scale=0.25]{REF_M7_Noisep1a.eps} & 
\includegraphics[scale=0.25]{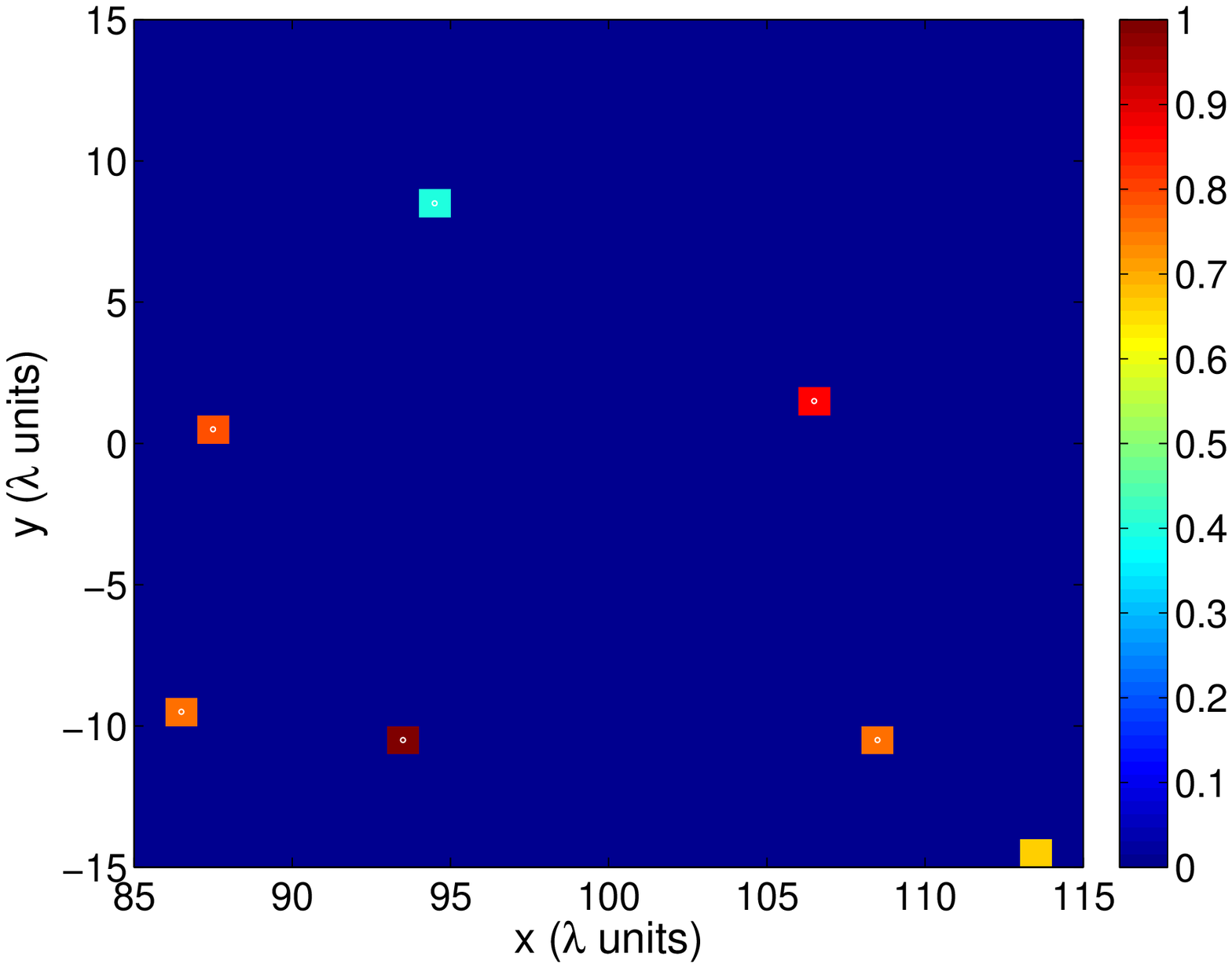} & 
\includegraphics[scale=0.25]{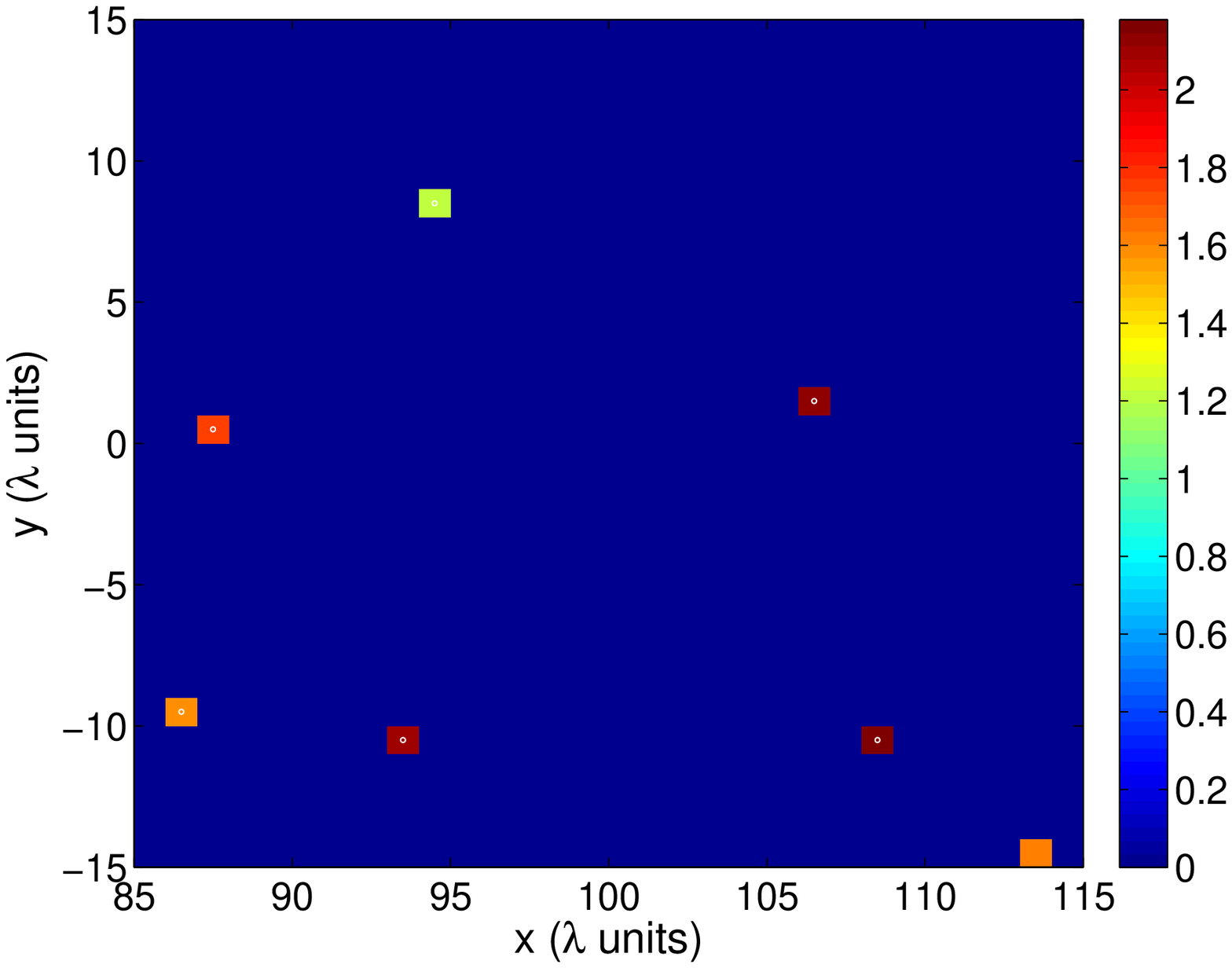}
\end{tabular}
\caption{$10 \%$ noise. The top and bottom rows show the images obtained with MUSIC and MMV, respectively. Top row from left to right (MUSIC):
original configuration of the scatterers, locations of the scatterers given by MUSIC, and amplitudes of the reflectivities obtained after nuclear norm minimization. Bottom row from left to right (MMV): original configuration of the scatterers, loctations of the effective sources,
and amplitudes of the reflectivities obtained after the second step \eqref{rhostar}.
}
\label{fig:10noise}
\end{figure}

Figure \ref{fig:10noise} illustrates the results with $10\%$ of noise added to the data. The left column displays the original configuration of the scatterers,
which is the same for both MUSIC (top row) and MMV (bottom row) reconstructions. In the top row,  the middle plot shows the
locations of the scatterers given by the MUSIC imaging function \eqref{MUSIC}. The right plot shows the final image, obtained once the reflectivities have been estimated by solving  the nuclear norm minimization problem \eqref{minimumnuclearnorm}. We observe very accurate scatterer locations, although oversmoothed in two of the scatterers. The bottom row displays the images obtained with the MMV formulation.
The  middle plot shows the
locations of the effective sources given by the solution to \eqref{MMV21noise}. The right image shows the final image obtained with MMV, once the reflectivities of the scatteters have been found in the second step. Both, the locations and the reflectivities of the scatterers 
obtained with the MMV formulation are very accurate.

\begin{figure}[t]
\centering
\begin{tabular}{ccc}
\includegraphics[scale=0.25]{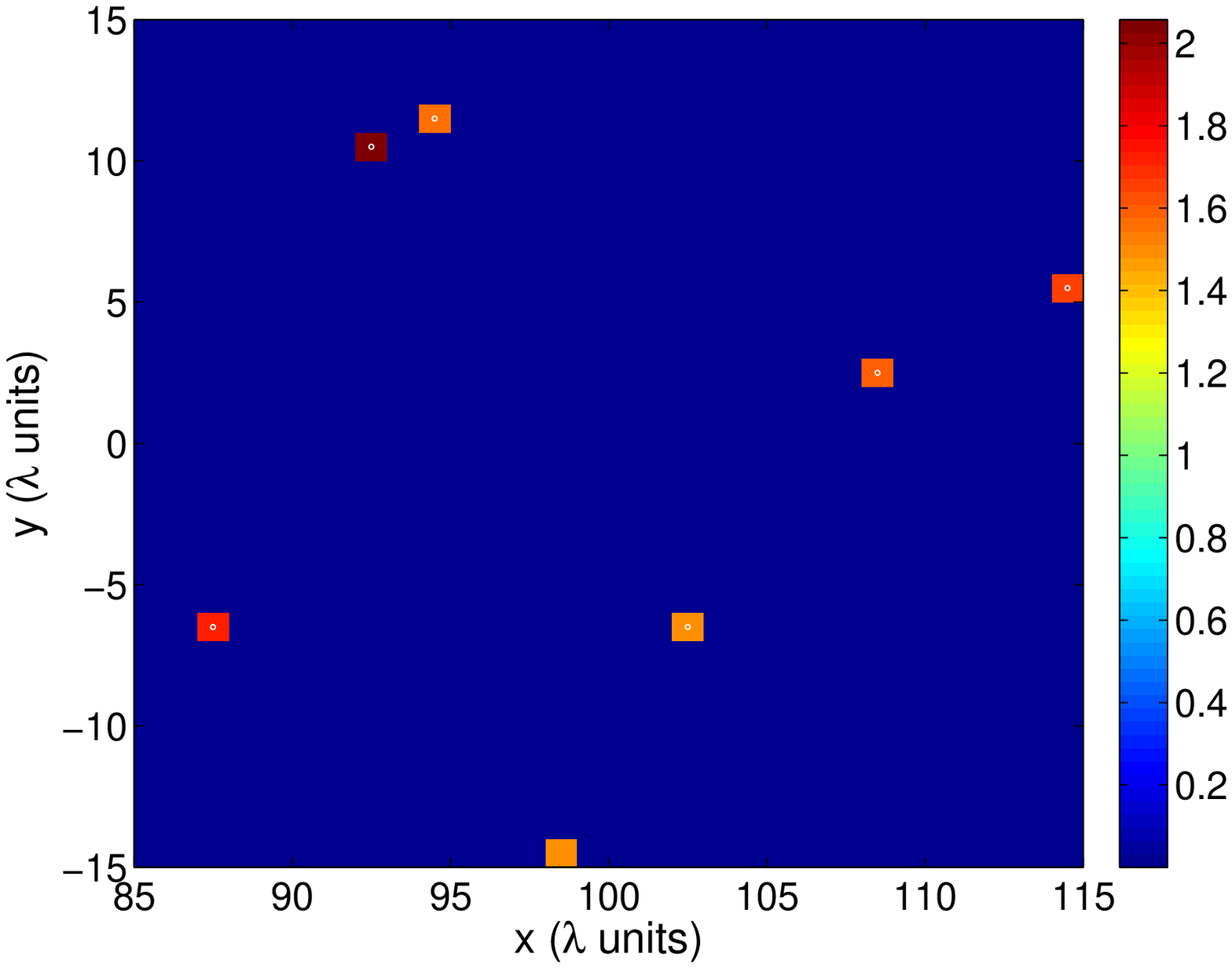} & 
\includegraphics[scale=0.25]{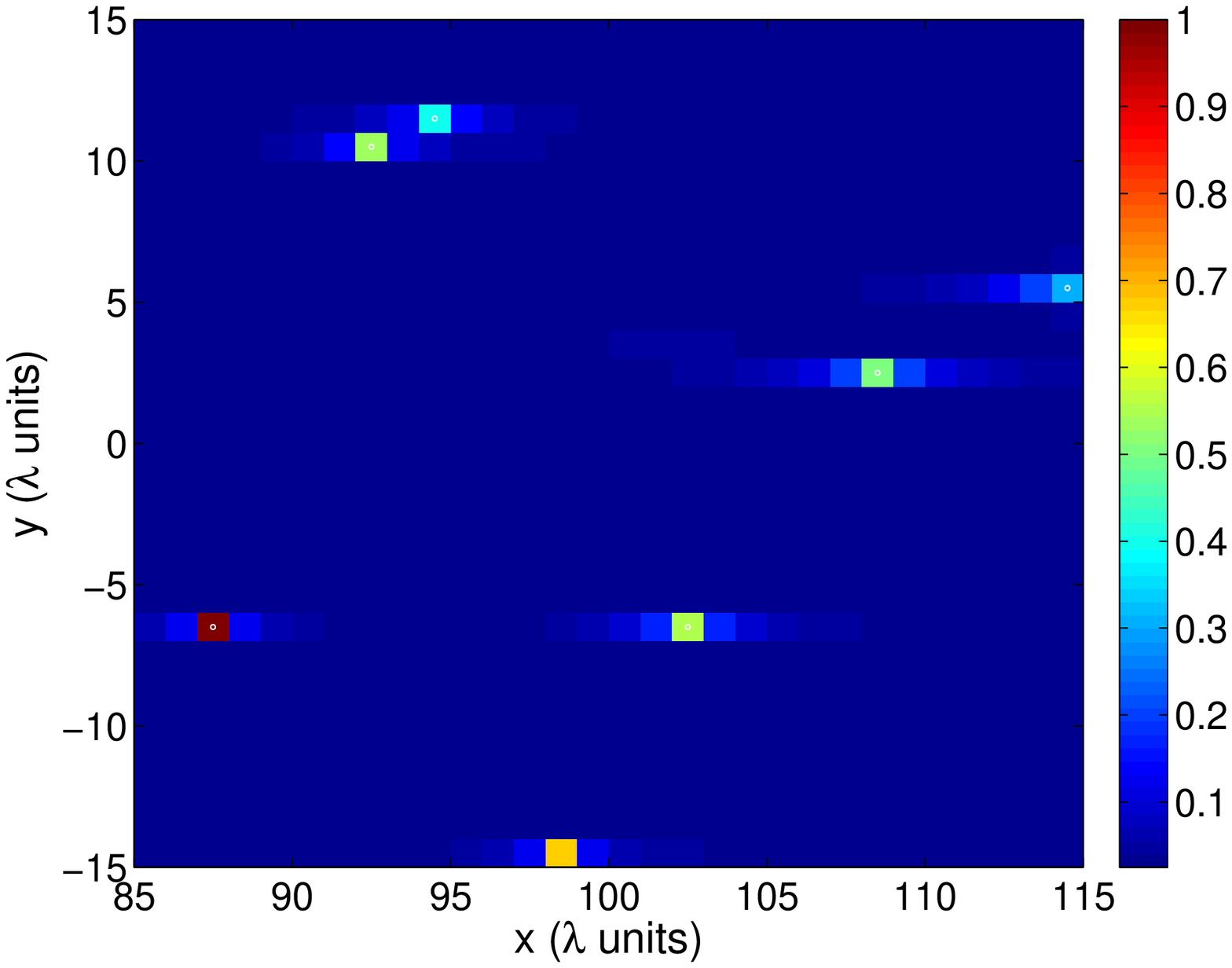} & 
\includegraphics[scale=0.25]{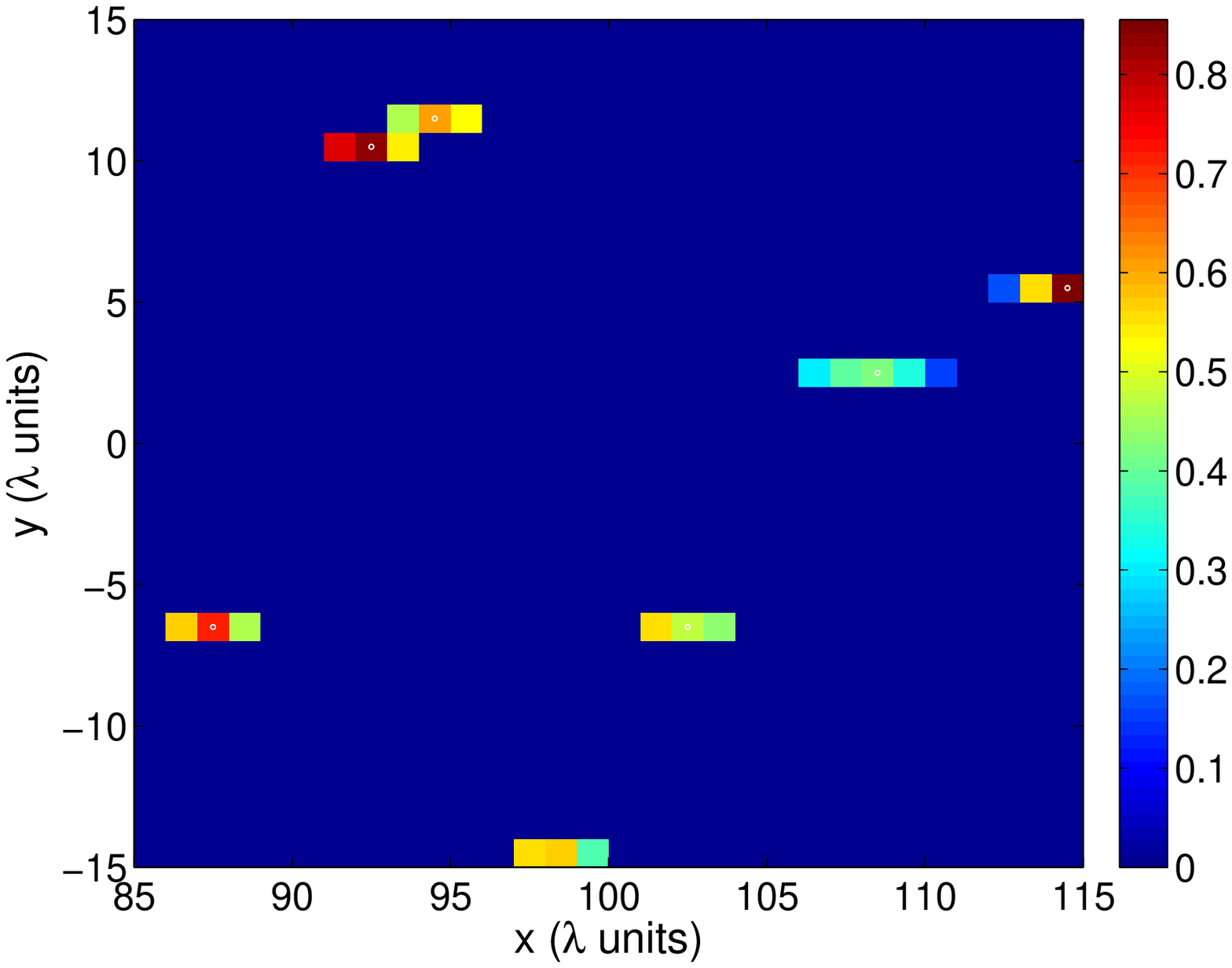} \\
\includegraphics[scale=0.25]{REF_M7_Noisep2a.eps} & 
\includegraphics[scale=0.25]{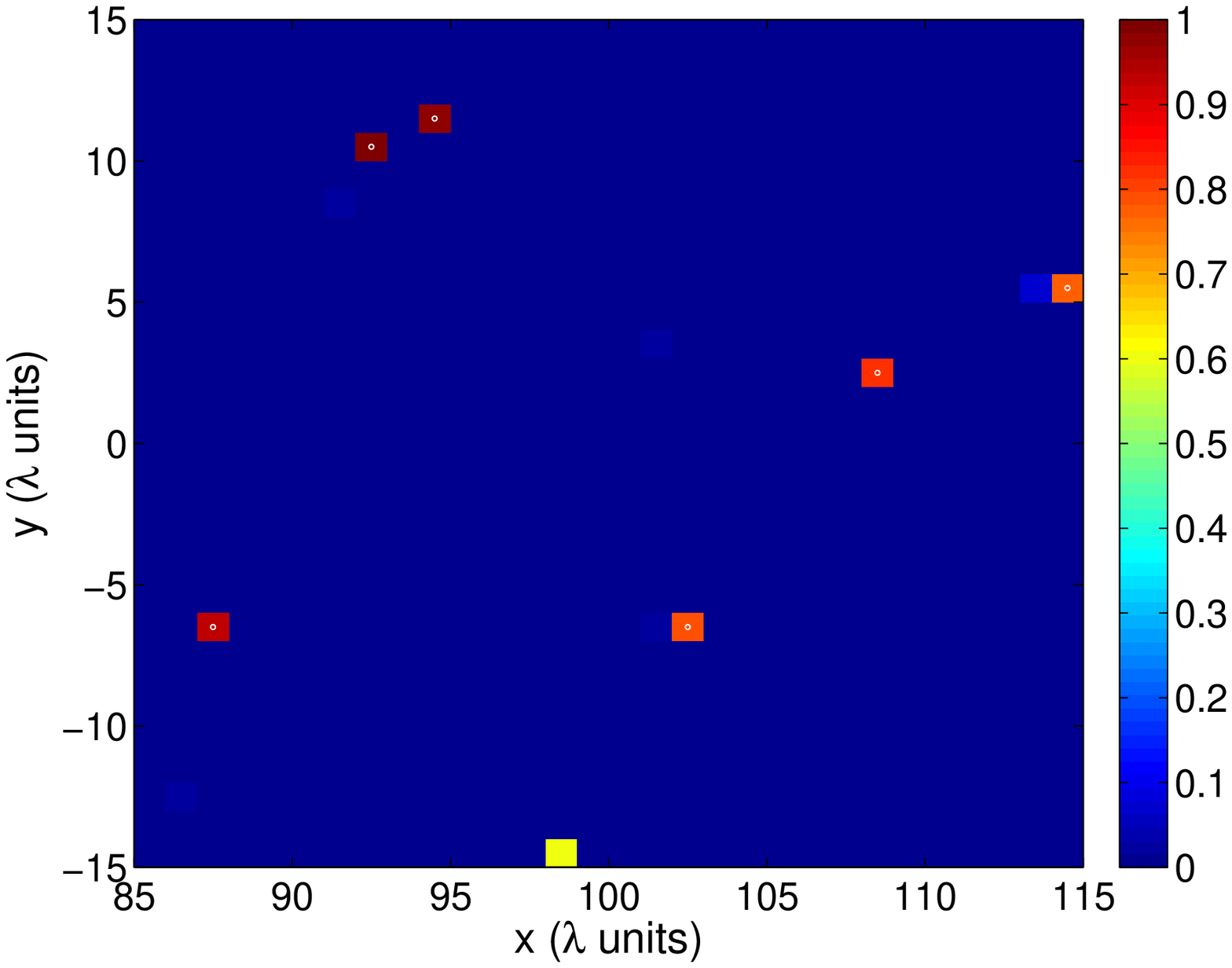} & 
\includegraphics[scale=0.25]{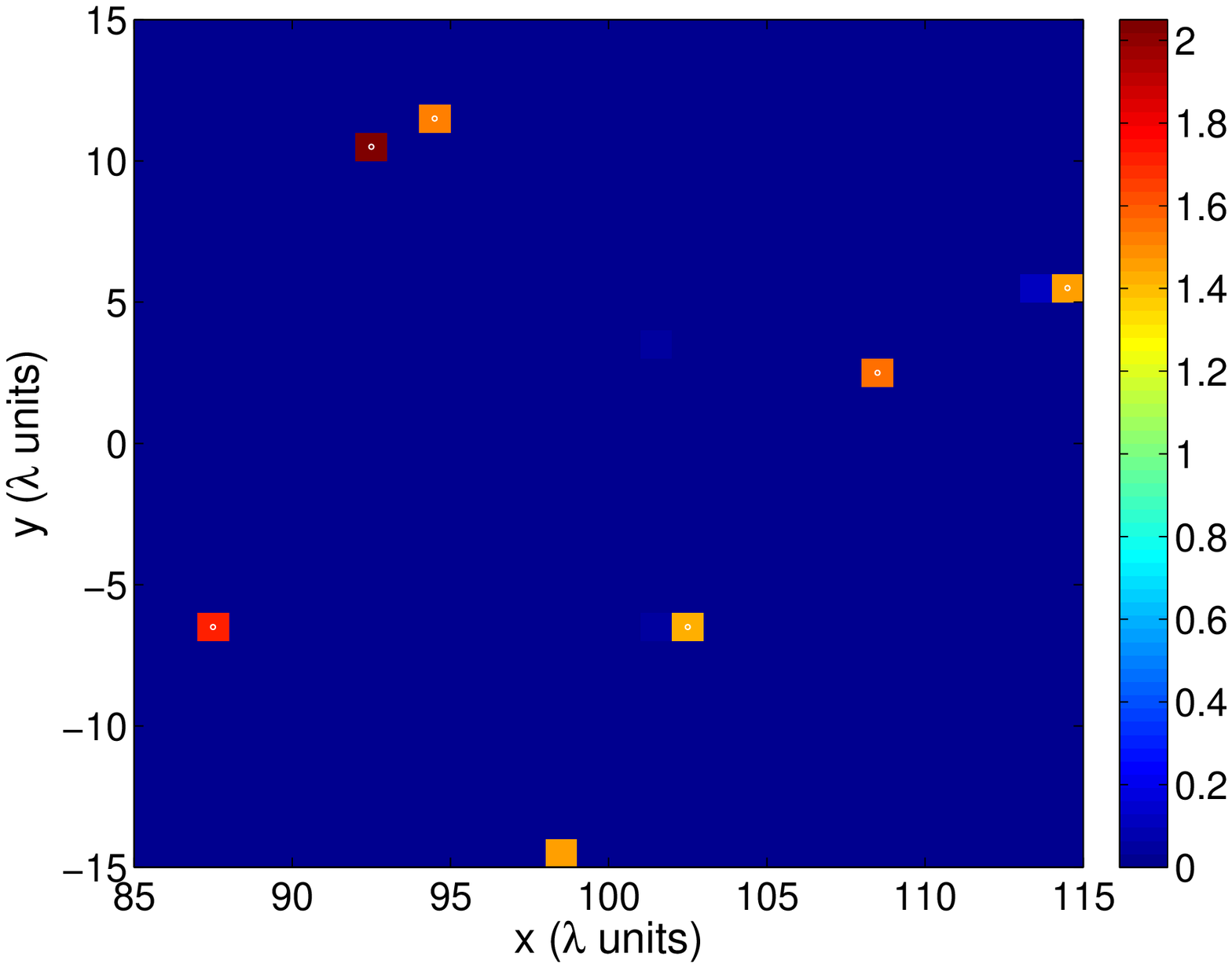}
\end{tabular}
\caption{
Same as Fig. \ref{fig:10noise} but with  $20 \%$ noise. 
}
\label{fig:20noise}
\end{figure}

Figure \ref{fig:20noise} is similar to Figure \ref{fig:10noise} but with $20\%$ of noise added to the data. The arrangement of
the images is the same as in that figure. The top row shows the results obtained with MUSIC, and
the bottom row the results obtained with MMV. Both methods still work well in locating the scatterers with $20\%$ of noise. The amplitudes of the reflectivities given by the MMV formulation are more accurate than those obtained
with MUSIC and nuclear norm minimization.

\begin{figure}[t]
\centering
\begin{tabular}{ccc}
\includegraphics[scale=0.25]{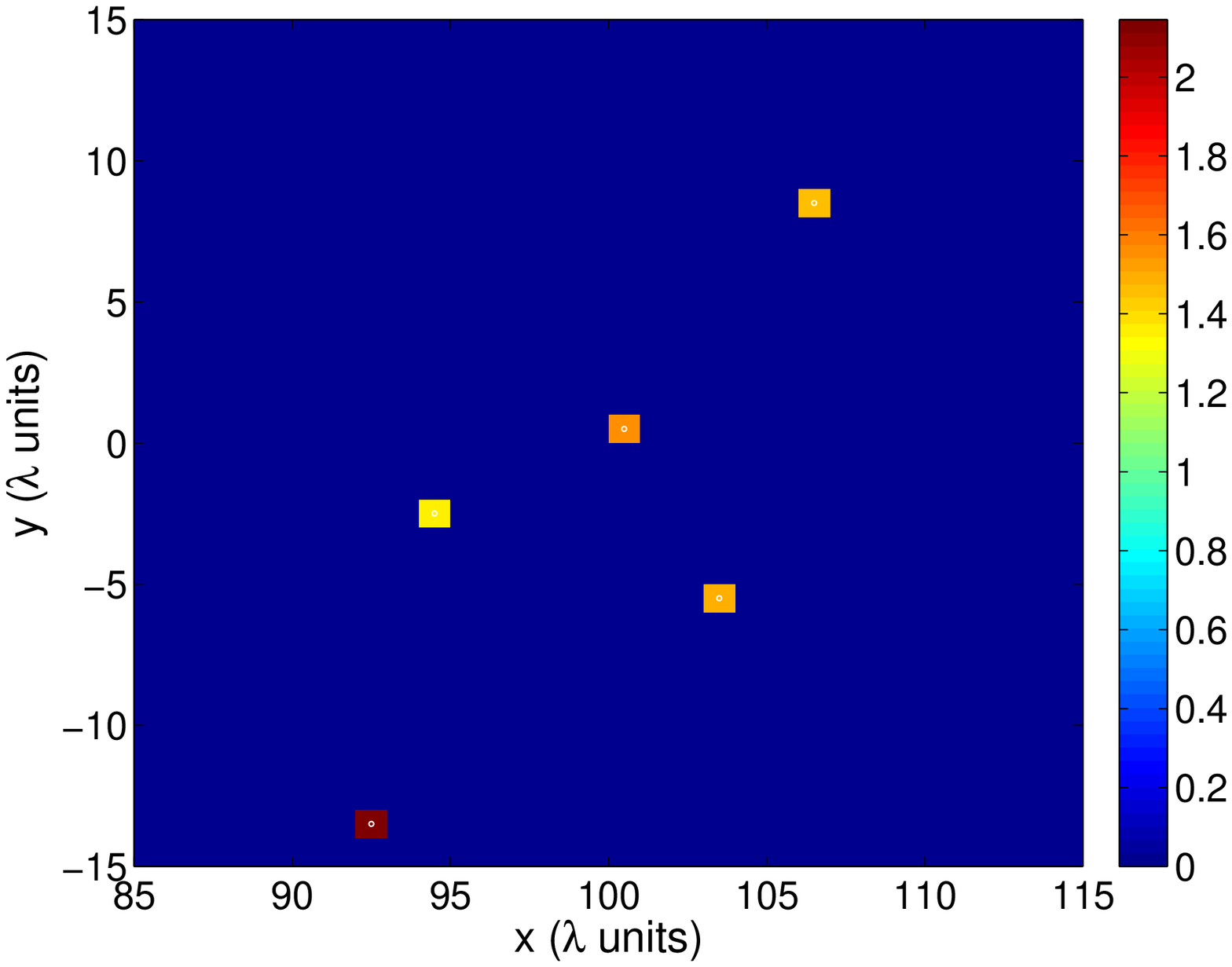} & 
\includegraphics[scale=0.25]{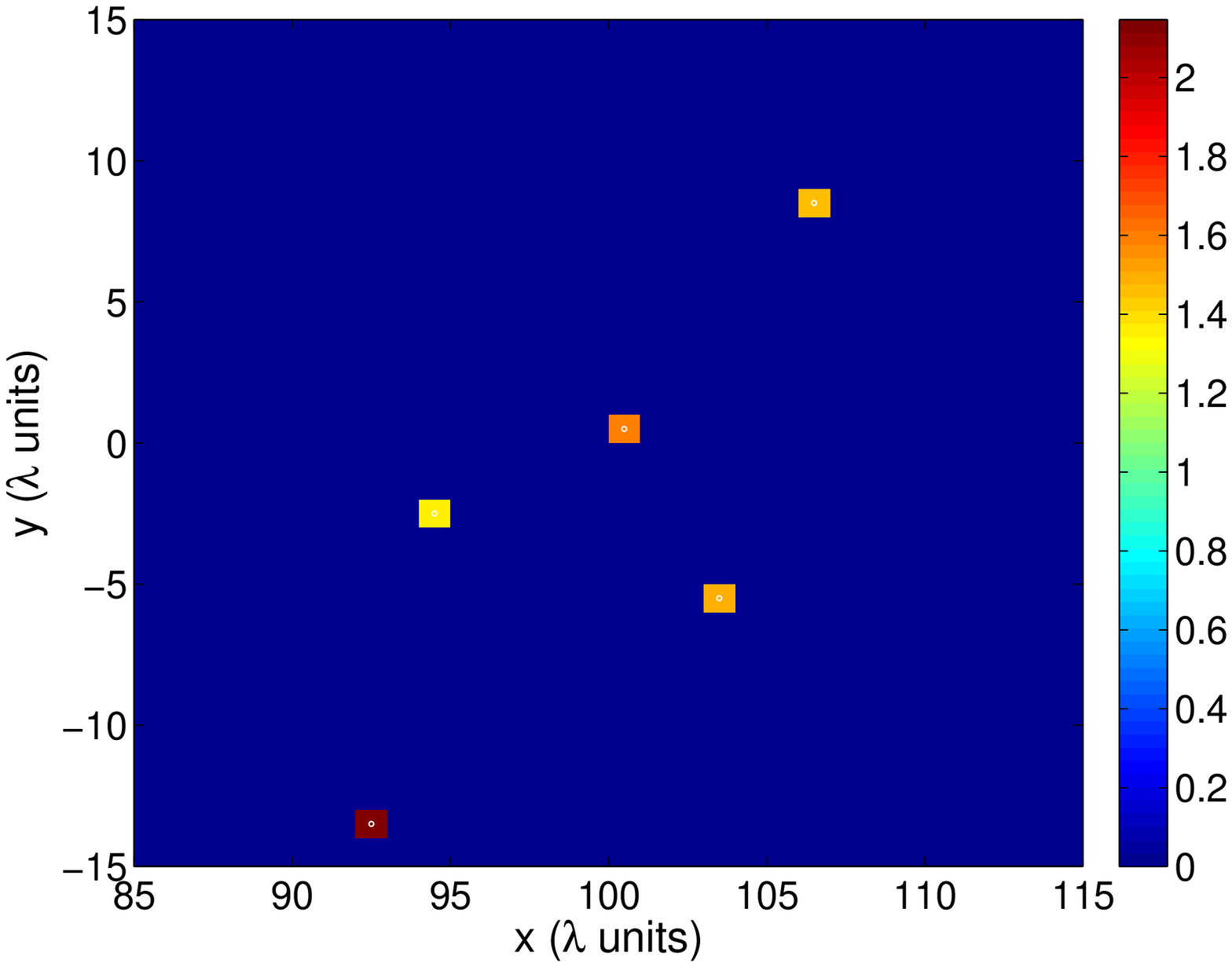} & 
\includegraphics[scale=0.25]{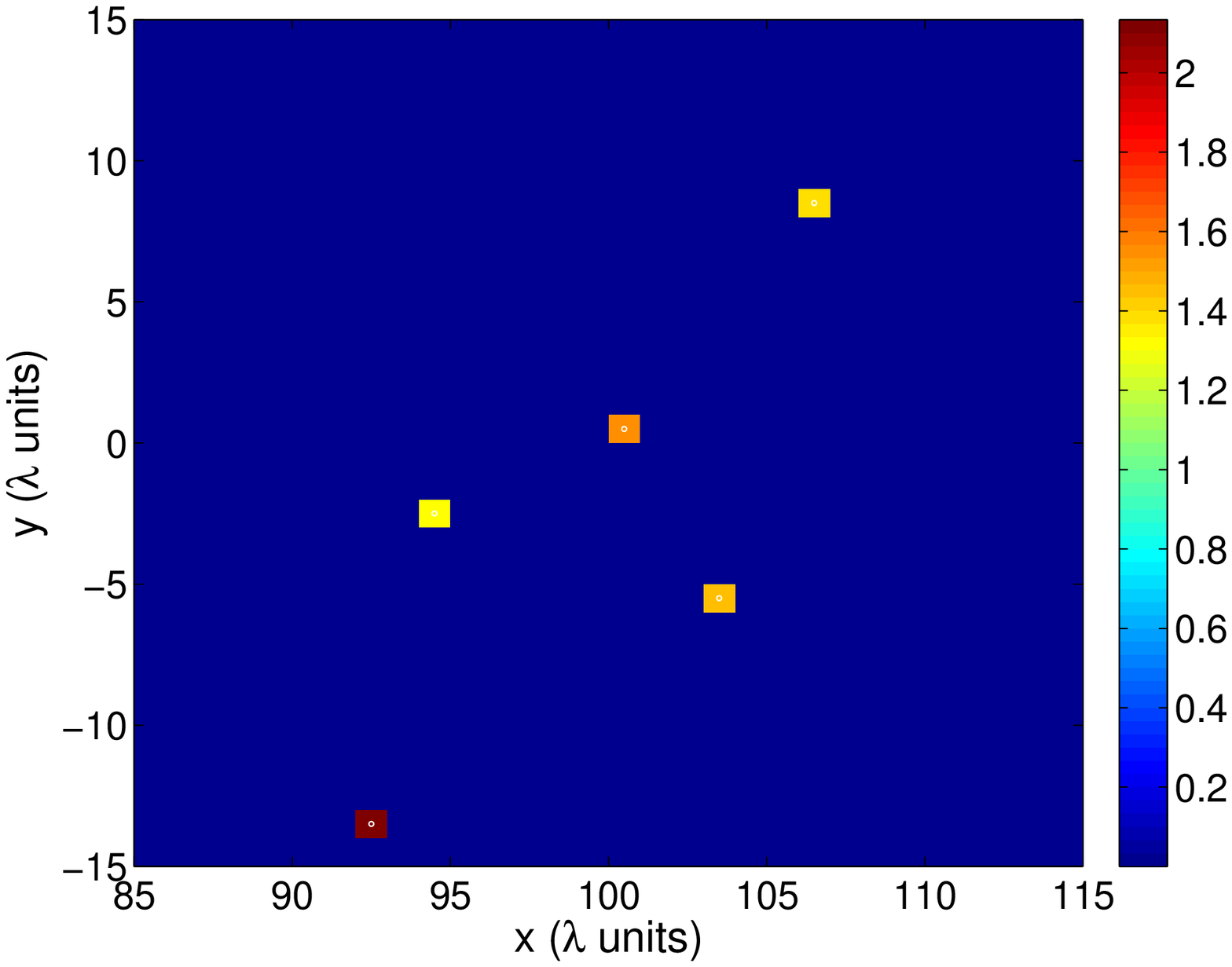} 
\end{tabular}
\caption{Incomplete set of illuminations. Only $50\%$ of the illuminations are used and $5 \%$ of noise is added to the data. The missing
entries of $\vect\wM(\omega)$ are found by matrix completion \eqref{minnuclear}.  The original configuration of the scatterers is shown in the left image, the image obtained with MUSIC in the middle image, and the imaged obtained with the MMV formulation in the right image. 
}
\label{fig:matrixcompletion}
\end{figure}

Next, we study the performance of the two methods with partial illumination, i.e., when the images are formed from an incomplete
set of illuminations as discussed in subsection \ref{sec:incompletedata}. First, we consider the case in which the data are corrupted.
Only data from some pairs of transducers, randomly selected, are available. In this case, the missing entries of $\vect\wM(\omega)$ are found by using matrix completion, i.e., by solving \eqref{minnuclear}.
Figure \ref{fig:matrixcompletion} shows the results when data from $50\%$ of the pairs of transducers in the array, randomly selected,  are used to form the images. This means that we have to recover the low rank data matrix $\vect\wM(\omega)$ from a random sampling of $50\%$ of its (noisy) entries. $5 \%$ of noise was added to the data in this experiment.
The distribution of scatterers to be recovered  is shown in the left plot, and the images obtained with MUSIC and MMV  in the middle and right plots, respectively. Both images are very good. 

We  note that, as expected, matrix completion does not work well with more than $50\%$ of the entries of $\vect\wM(\omega)$ missing,
even with noiseless data. 
This is in agreement with the theoretical results on the number of randomly sampled entries required to reconstruct an unknown low rank matrix \cite{Candes12}.

\begin{figure}[t]
\centering
\begin{tabular}{cc}
\includegraphics[scale=0.25]{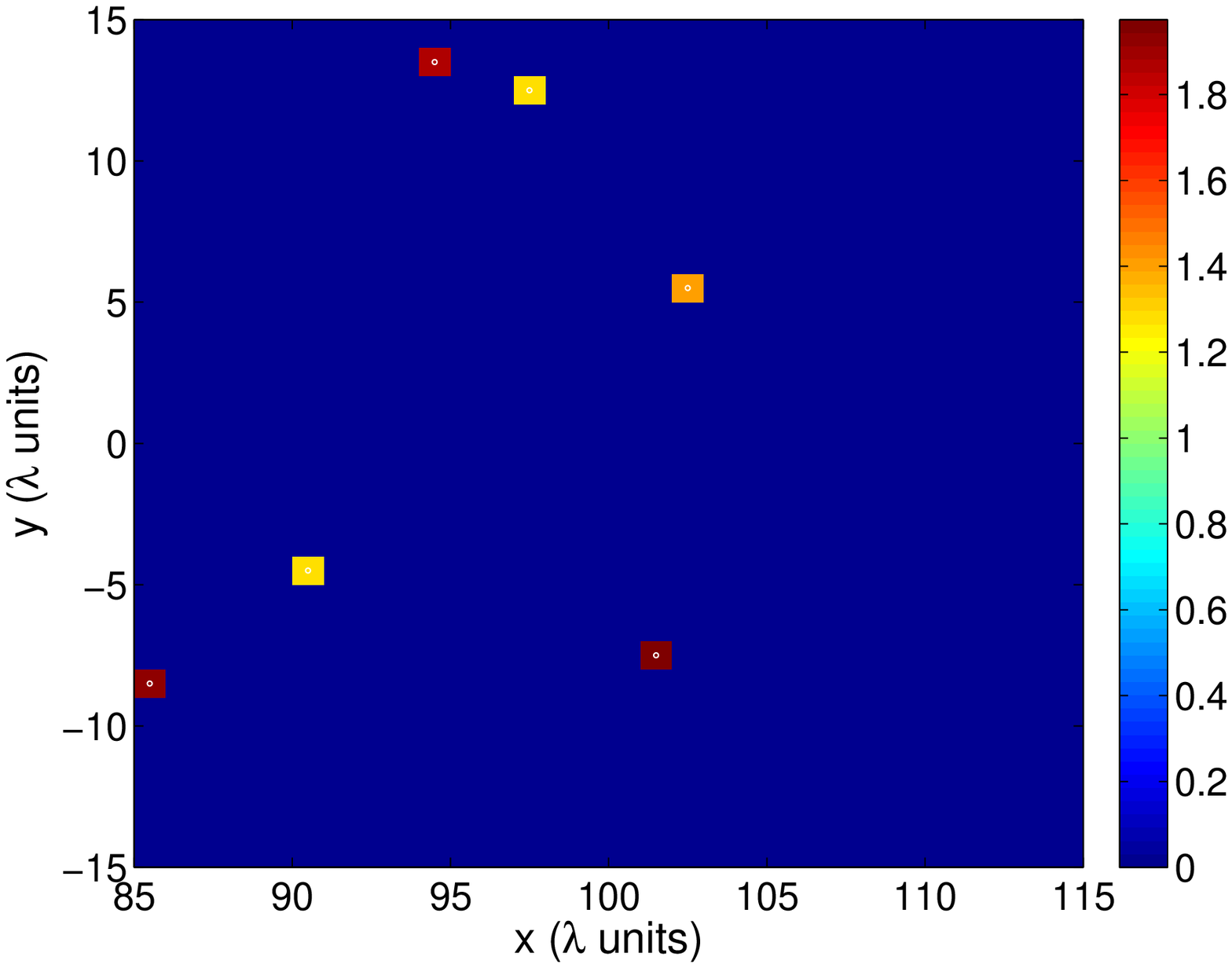} &
\includegraphics[scale=0.25]{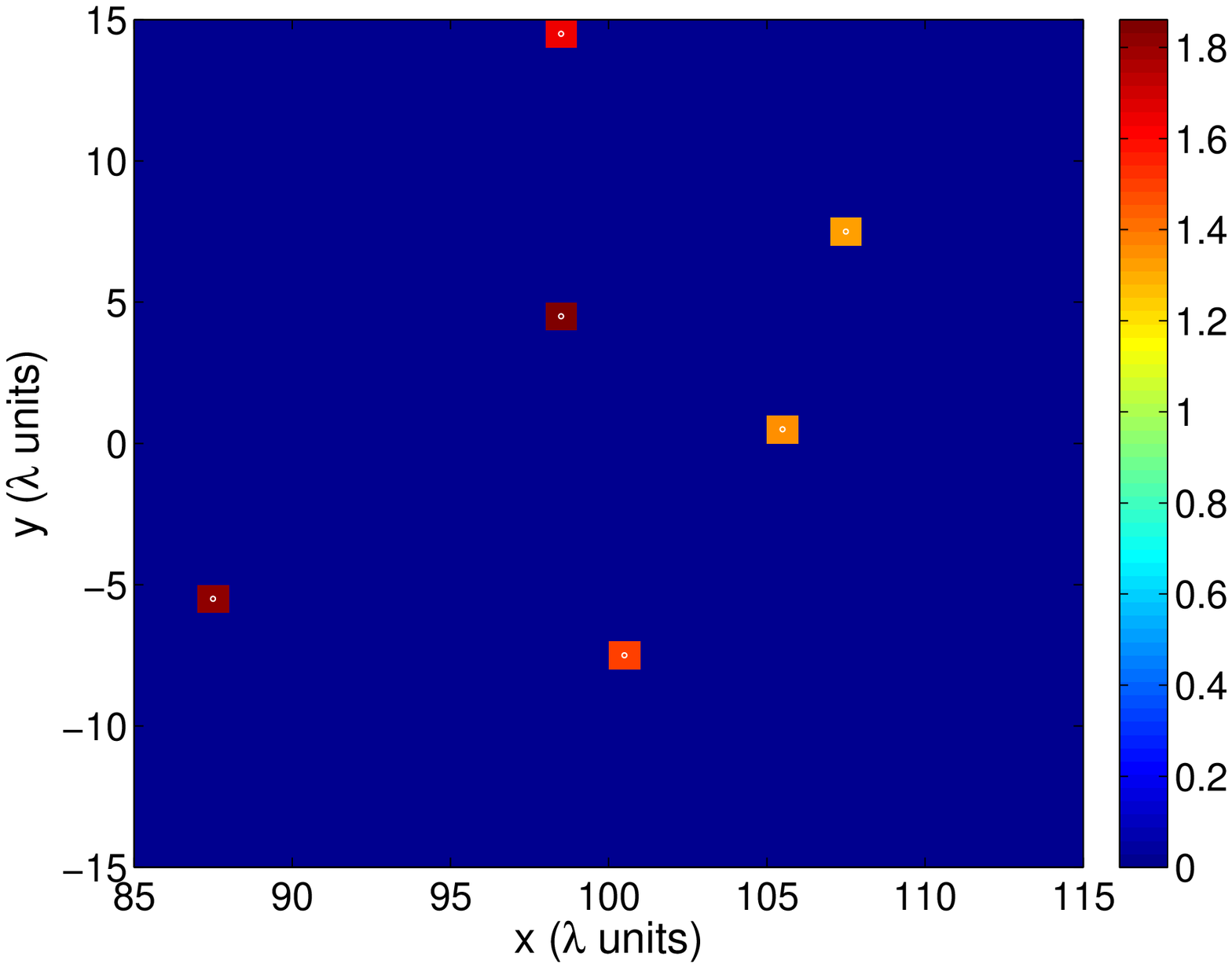}
\end{tabular}
\caption{Original configurations of the scatterers used for the numerical experiments shown in Fig. \ref{fig:edges0noise} (left image), and Figs.  \ref{fig:edges5noise}, and \ref{fig:edges10and20noise} (right image).
}
\label{fig:ref_for_edges}
\end{figure}

Finally, we examine the results when only a few transducers at the edges of the array are used to illuminate the IW. 
In this case,  intensity-only imaging is applied directly to the matrix formed by the submatrices at the four corners of $\vect\wM(\omega)$, without matrix completion.  Figure \ref{fig:ref_for_edges} shows two reference images used for the  study of the performance of MUSIC and MMV when illumination from the edges of the array is used.  In the next experiments, we only show the location of the scatterers recovered by these two methods. We do not carry out the second steps to estimate the reflectivities of the scatterers.

\begin{figure}[t]
\centering
\begin{tabular}{ccc}
\includegraphics[scale=0.25]{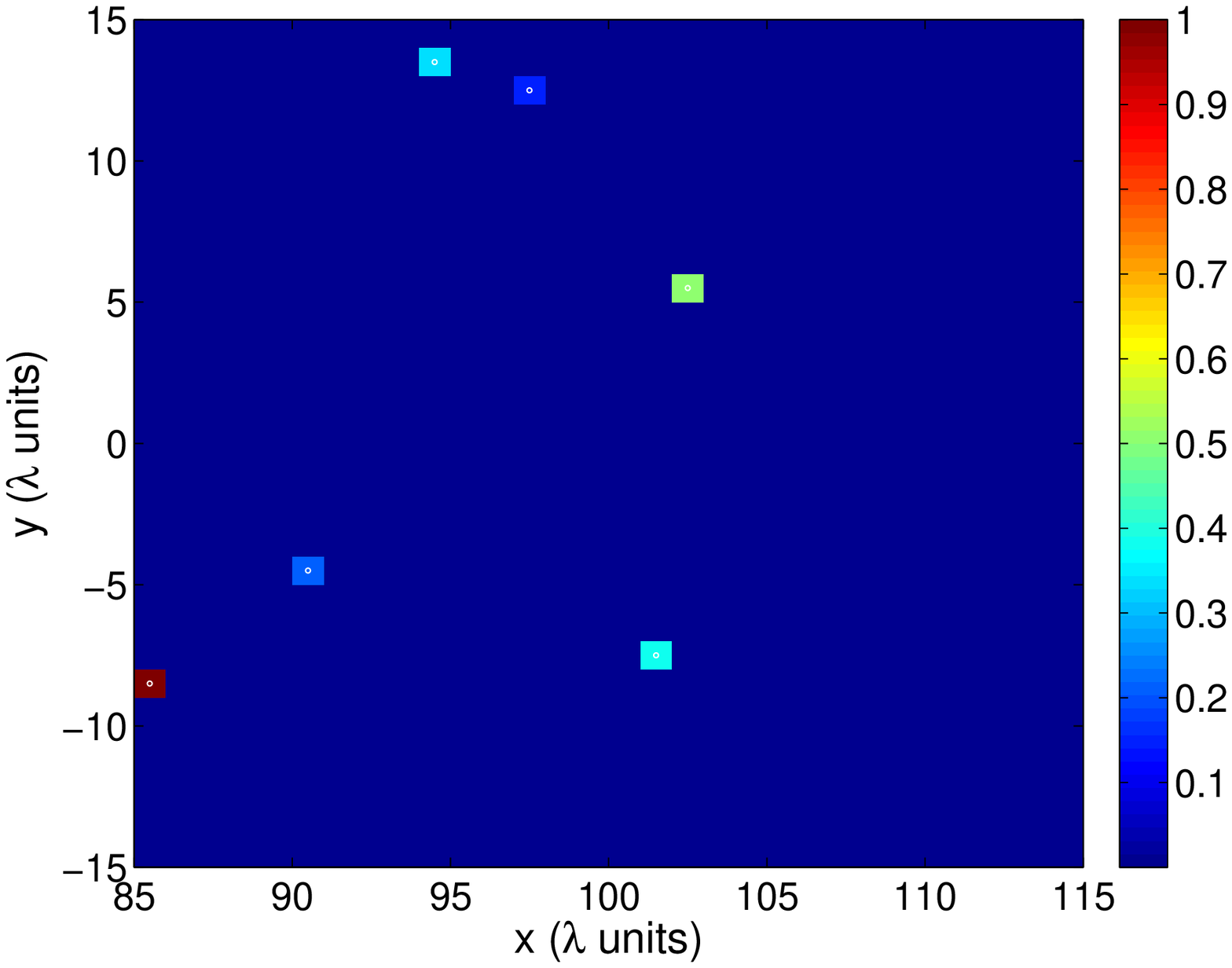} & 
\includegraphics[scale=0.25]{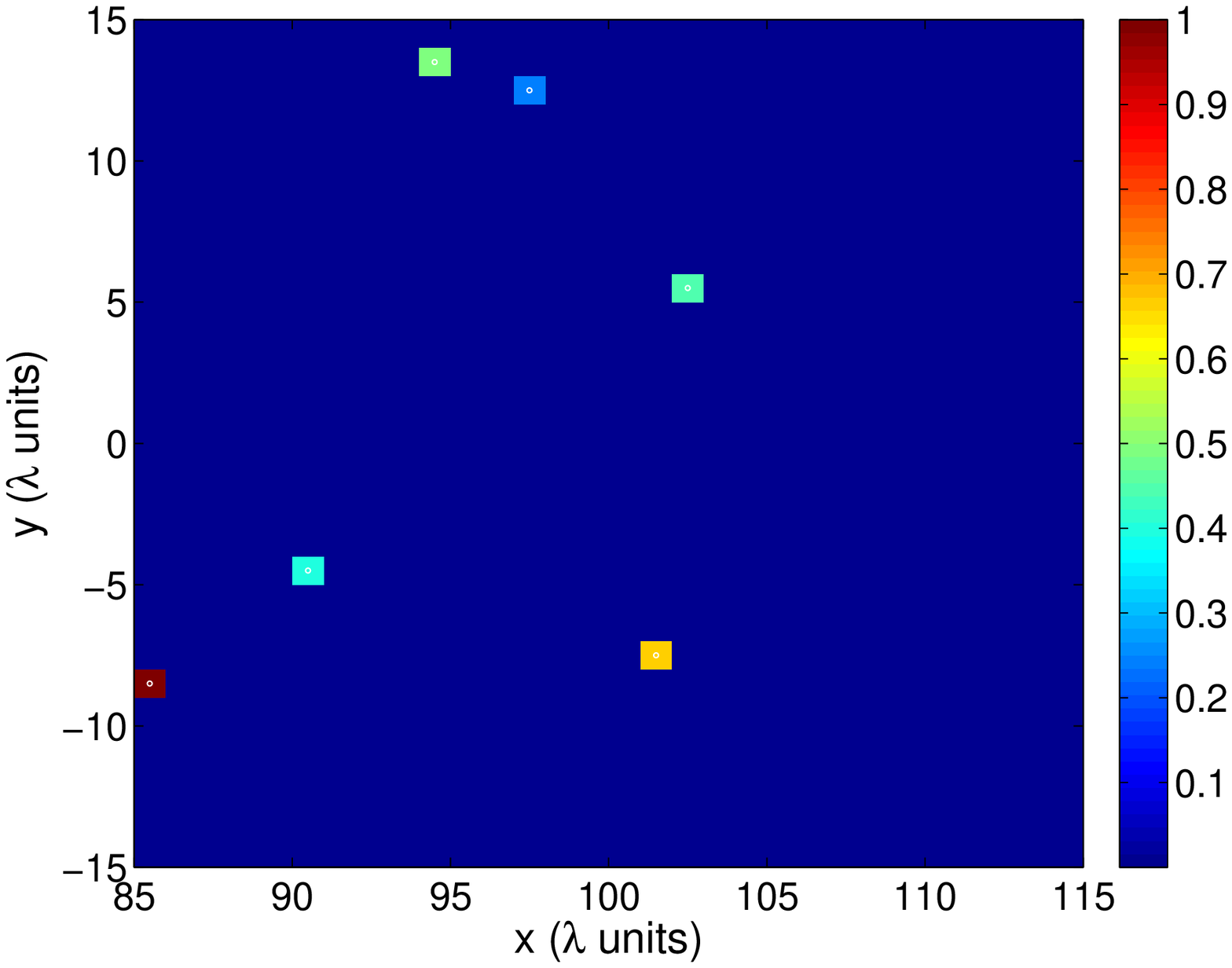} & 
\includegraphics[scale=0.25]{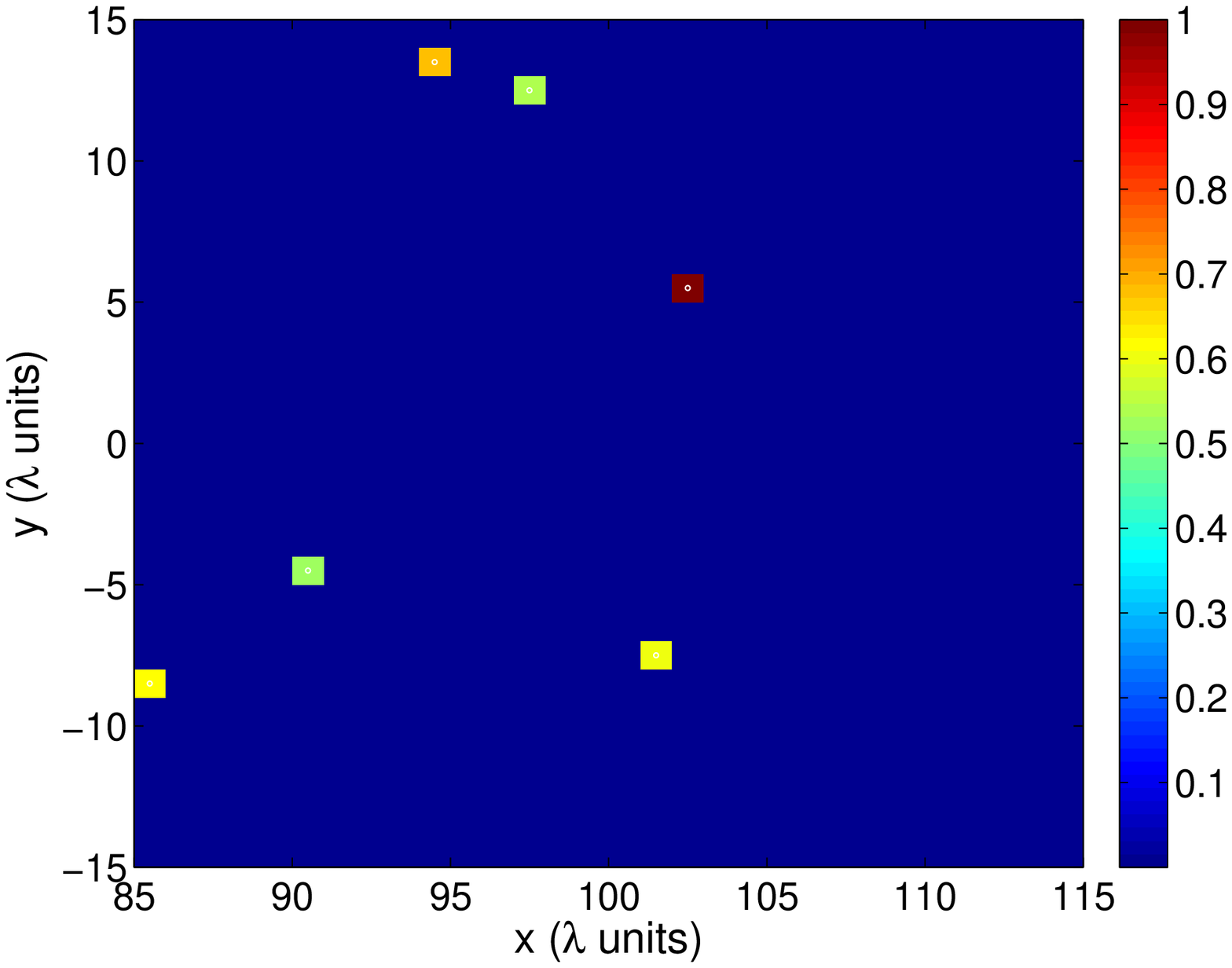} \\
\includegraphics[scale=0.25]{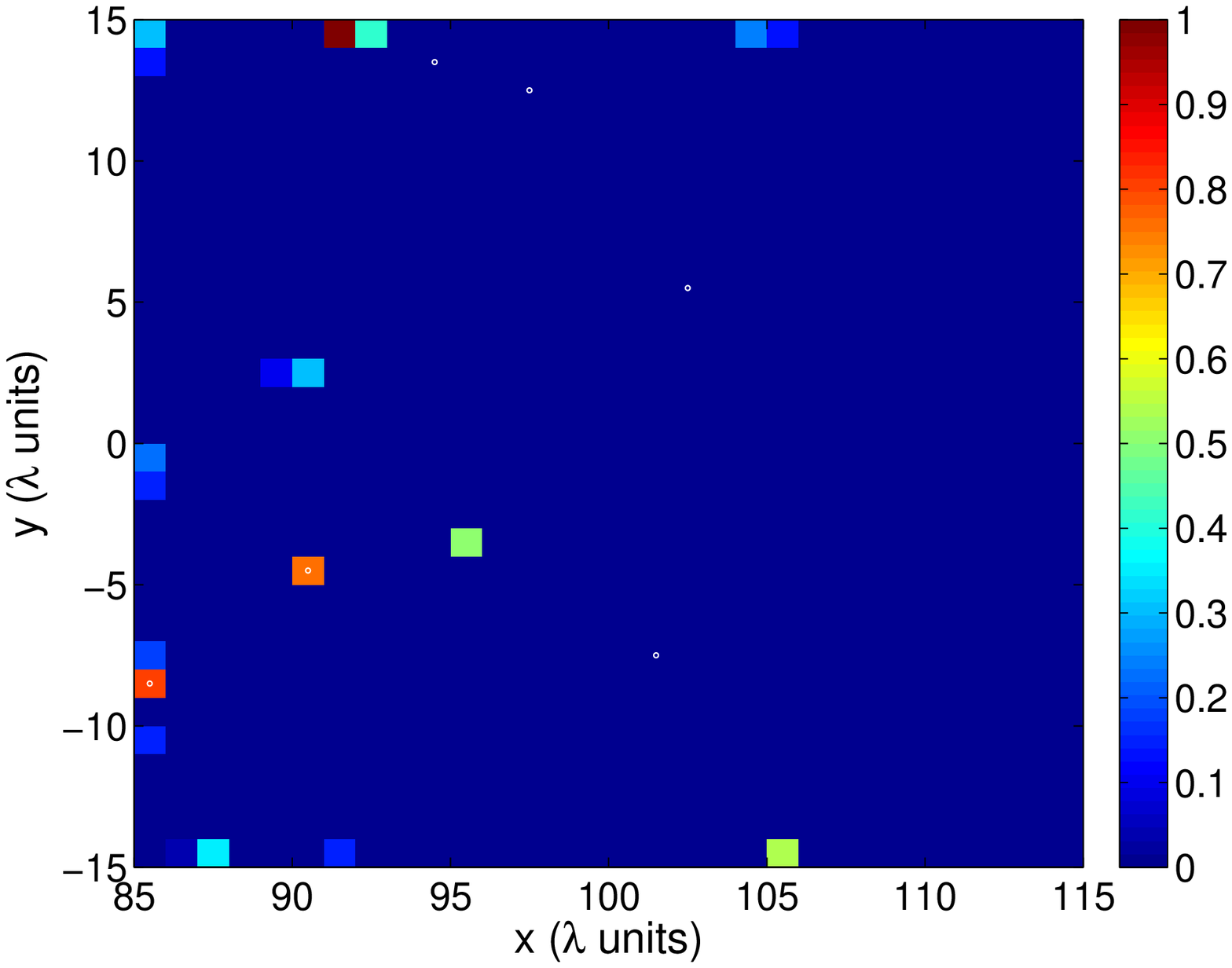} & 
\includegraphics[scale=0.25]{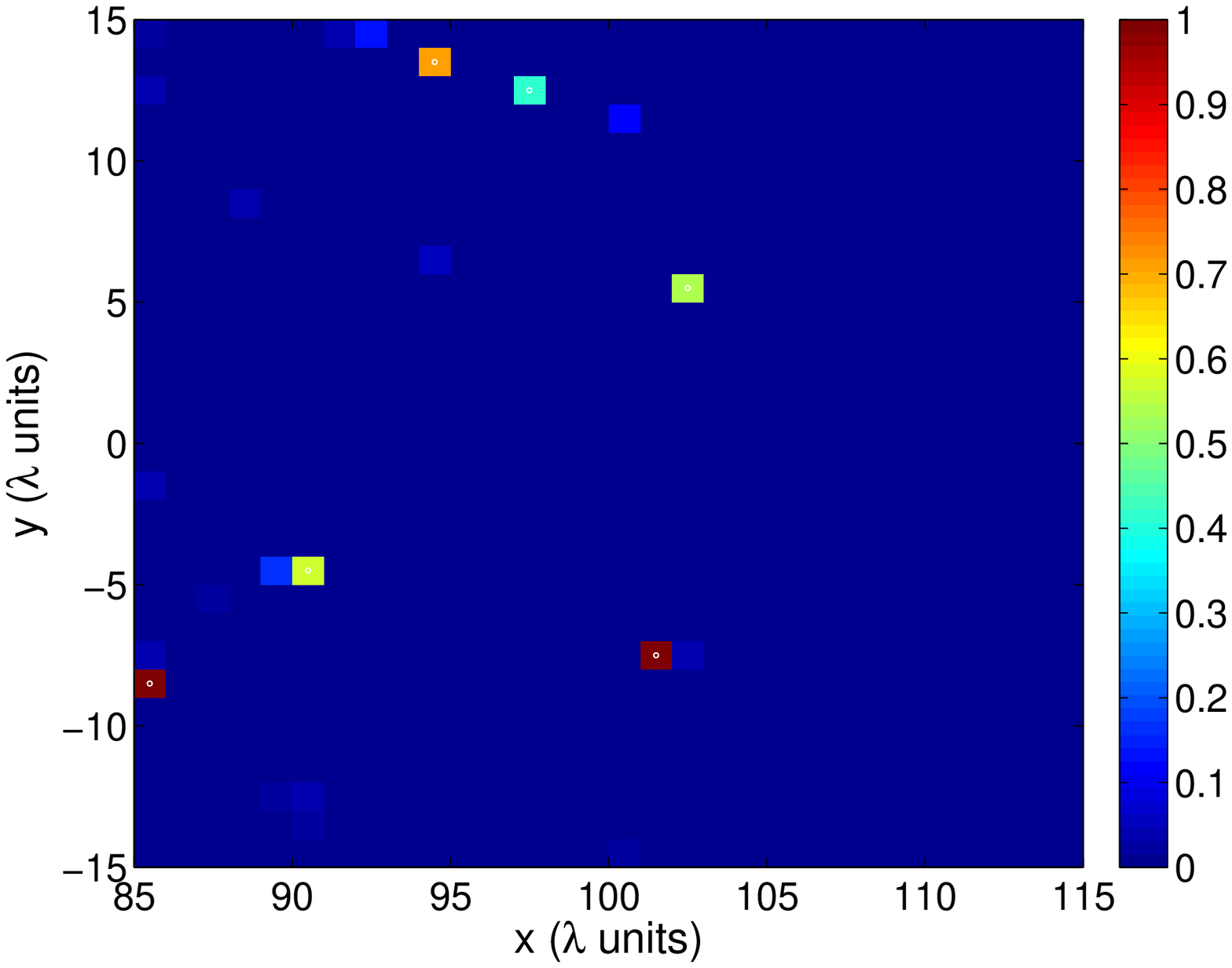} & 
\includegraphics[scale=0.25]{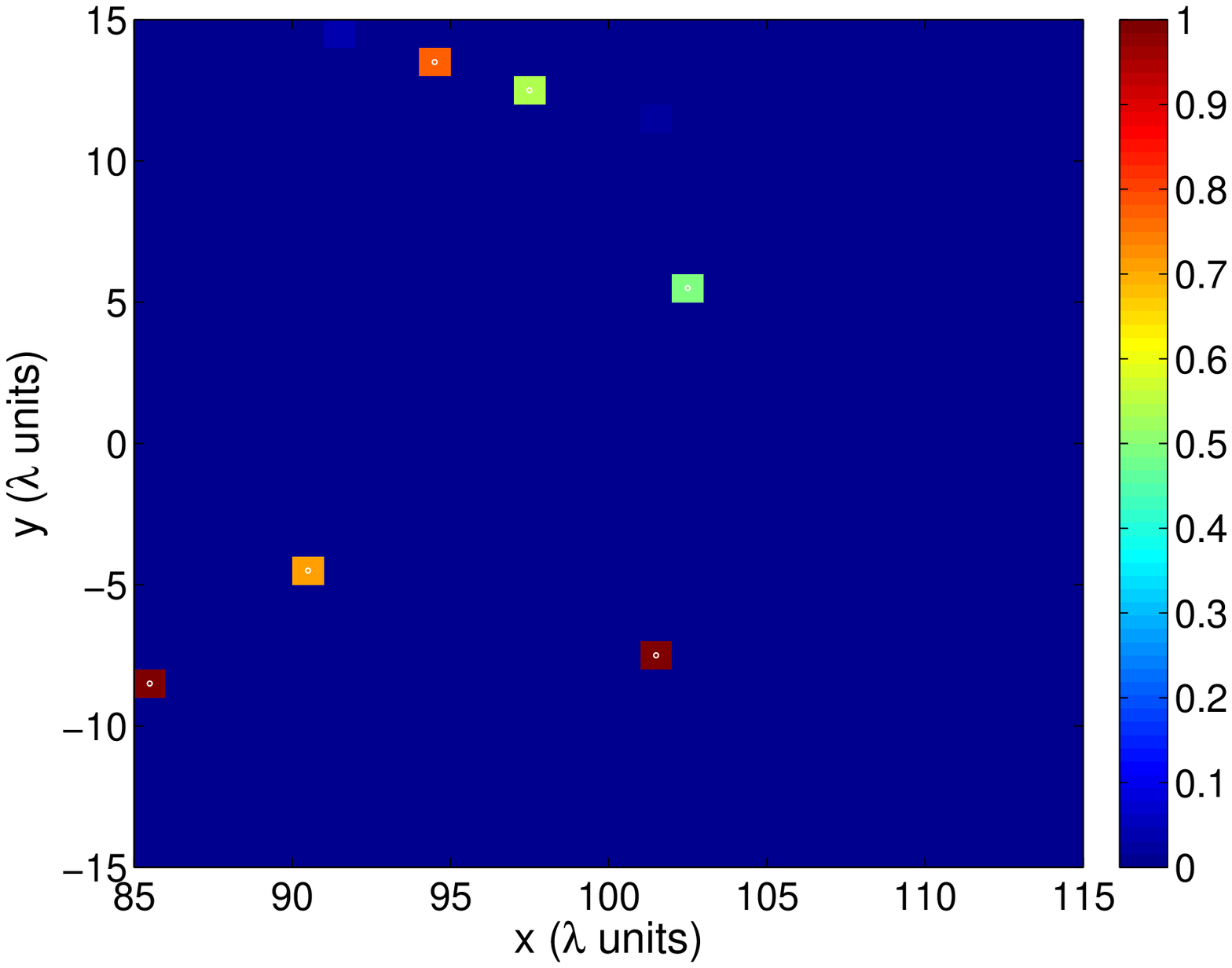}
\end{tabular}
\caption{Incomplete set of illuminations with no noise in data. Only partial illumination from the edges of the array is used.
$4$ (left column), $16$ (middle column), and $28$ (right column) transducers at each edge of the array are used.
The original configuration of the scatterers is shown in the left image of Fig. \ref{fig:ref_for_edges}.
Top row: location of the scatterers obtained with MUSIC. Bottom row:  location of the scatterers obtained with MMV. 
}
\label{fig:edges0noise}
\end{figure}

In the top row of Figure \ref{fig:edges0noise} we show the locations of the scatterers given by MUSIC when 
$4$ (left image), $16$ (middle image), and $28$ (right image) transducers at each edge of the array are active and illuminate the image window. There is no noise in data in these experiments. The original configuration of the scatterers is displayed in the left image of Figure  \ref{fig:ref_for_edges}.  
It is remarkable that only a few transducers at the edges of the array are enough to find the location of the scatterers accurately using 
MUSIC when there is no noise in the data. In fact, even with only $N_{active}=8$ transducers ($4$ at each edge of the array) MUSIC 
locates the scatterers accurately. This is so because the image is sparse, with only $M=6$ scatterers in the image window, and 
$N_{active}>M$ transducers are enough to compute the signal and noise subspaces, where $N_{active}$ is the number of 
transducers used during the illumination process. We note, though, that  the peaks  are sharper at all 
the scatterer locations when more transducers are used. Hence, it is expected that the robustness of MUSIC with respect to  noise 
increases when more transducers are used.

In the bottom row of Figure \ref{fig:edges0noise}, we show the locations of the scatterers given by the MMV approach. It is
apparent that the MMV approach is not able to find the locations of the scatterers using only a few transducers. More data are 
necessary to achive good results with MMV. We remind that, through the polarization identity, the MMV approach uses complete data, 
including phases, only at those (pairs) of transducers used to illuminate the image window. Hence, the less pairs of transducers are used, 
the less data are available for MMV and the less constrains there are in \eqref{MMV21noise}. 
Indeed,  the bottom left image in Figure \ref{fig:edges0noise} shows that MMV completely fails to locate the scatterers using 
$4$ transducers at the each edge of the array, and the bottom 
middle image shows a few ghosts using $16$ transducers, even though there is no noise in the data. Only with $28$ transducers, 
around $50 \%$ of the transducers in the array, the image obtained with MMV is accurate (right image in Figure \ref{fig:edges0noise}). 
Hence, we observe that when only a few transducers at the edges of the array are used to illuminate the IW, MUSIC is the preferred 
method for intensity-only imaging.

\begin{figure}[t]
\centering
\begin{tabular}{ccc}
\includegraphics[scale=0.25]{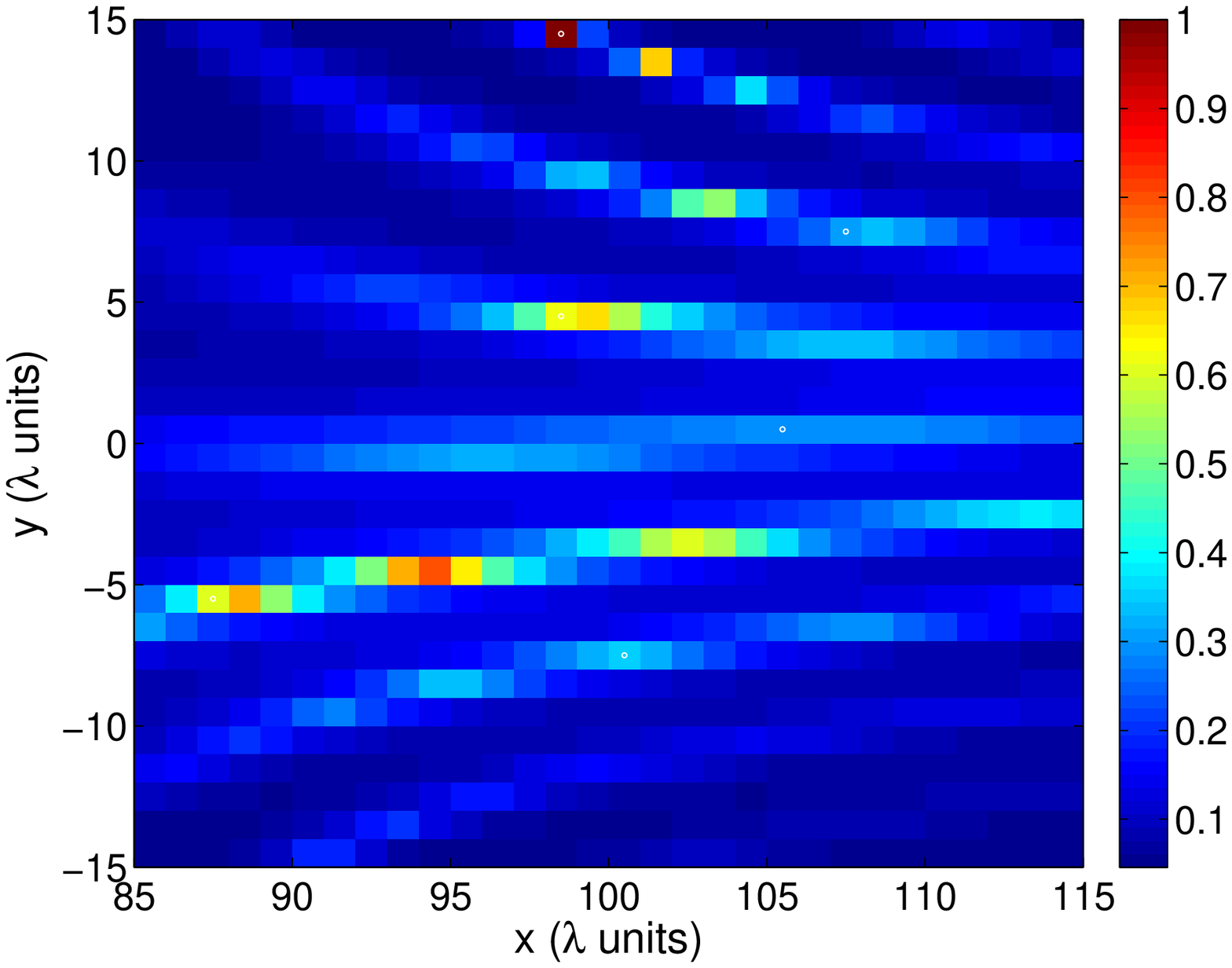} & 
\includegraphics[scale=0.25]{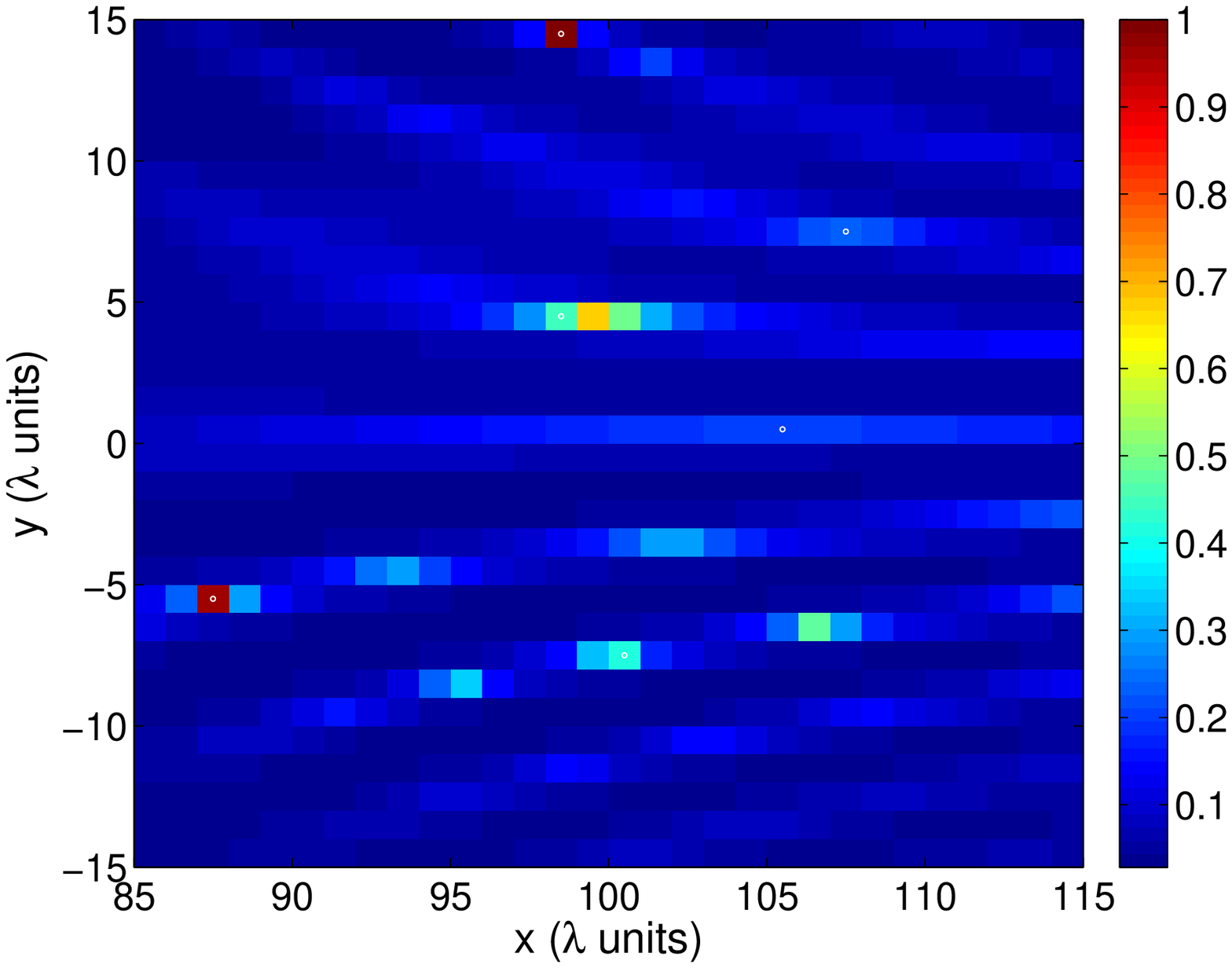} & 
\includegraphics[scale=0.25]{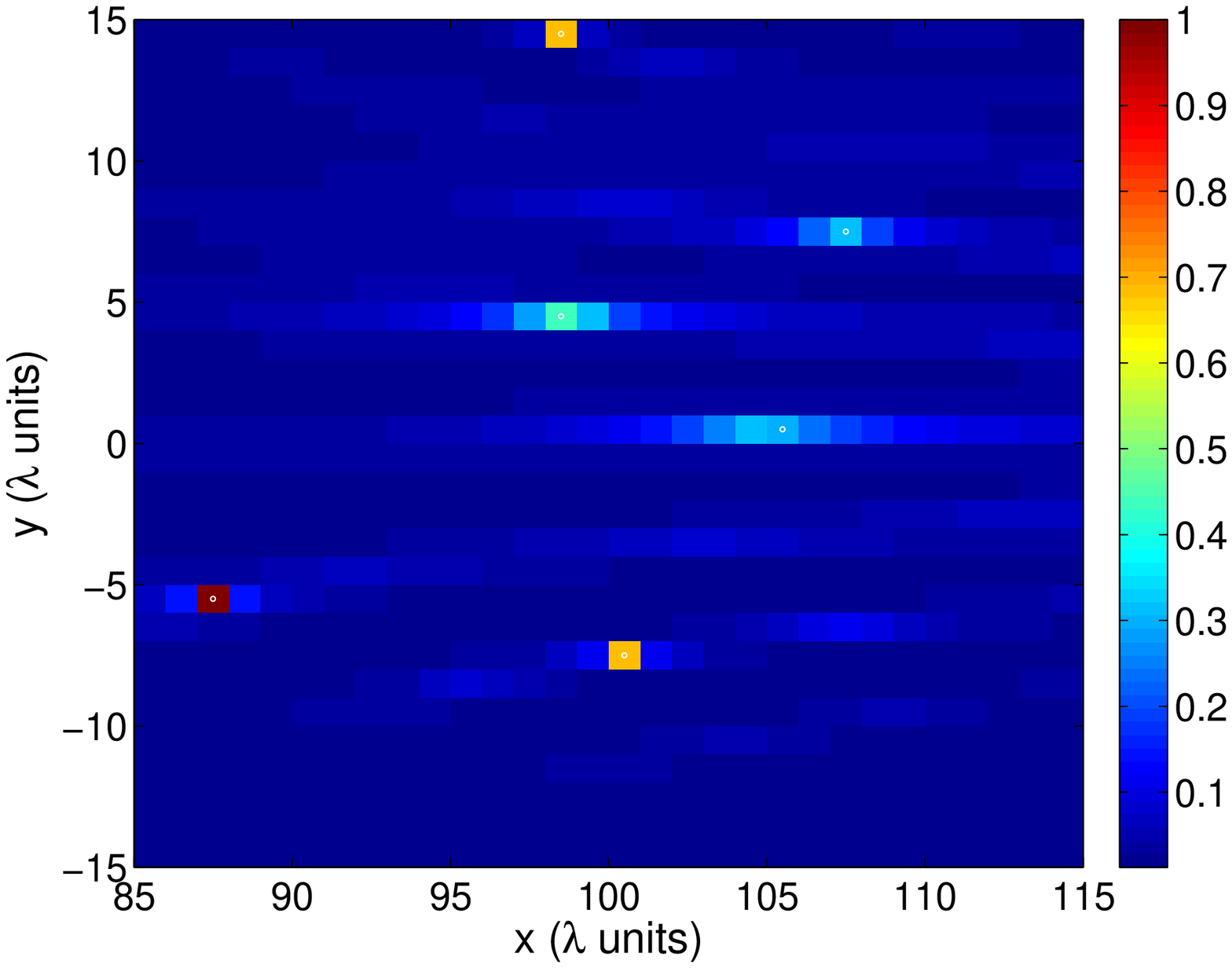} \\
\includegraphics[scale=0.25]{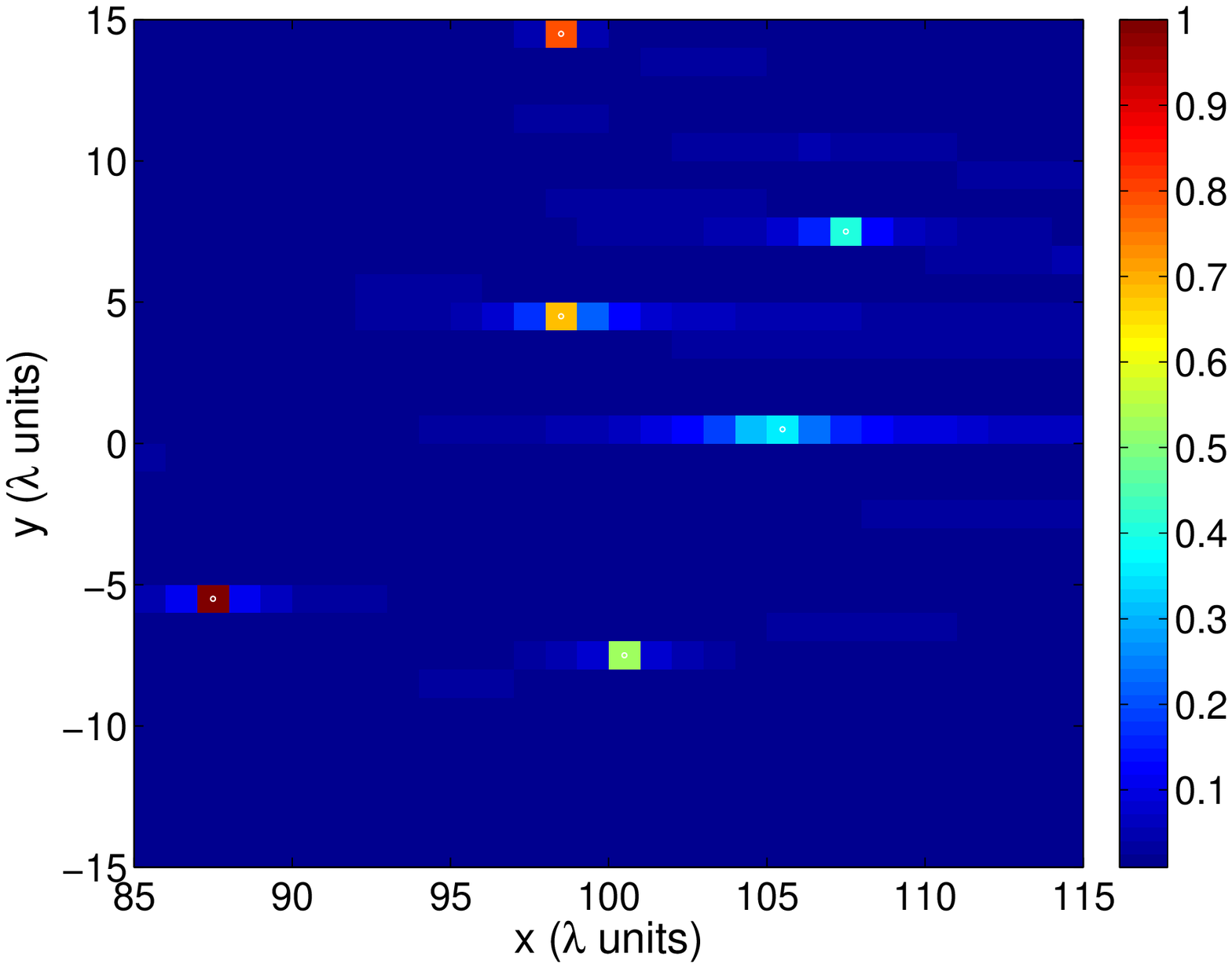} & 
\includegraphics[scale=0.25]{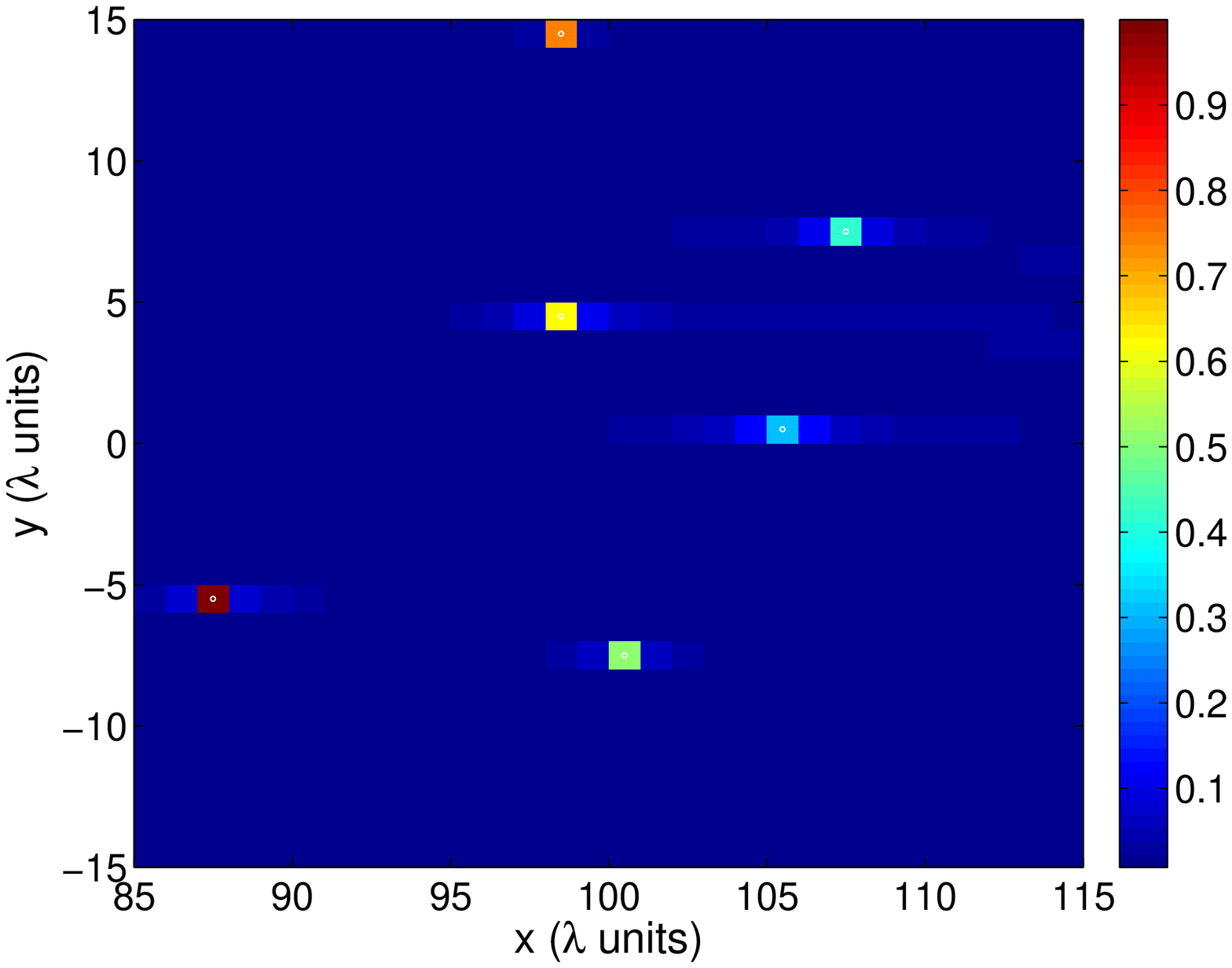} & 
\includegraphics[scale=0.25]{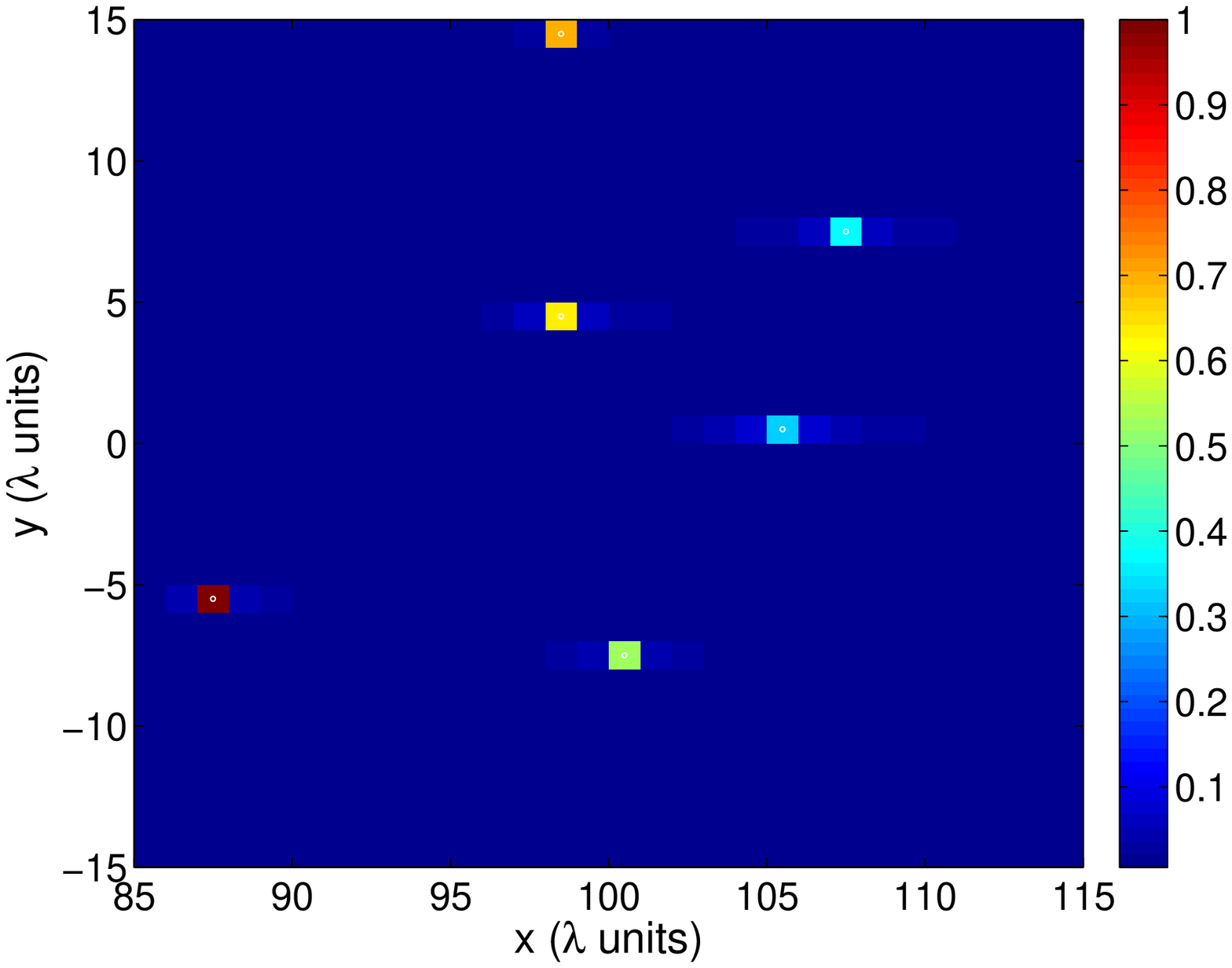}
\end{tabular}
\caption{Partial illumination from the edges of the array
($5 \%$ of noise is added to the data). From left to right and from top to bottom $4$, $8$, $12$, $16$, $20$ and $24$ transducers at each edge of the array illuminate the image window.  The locations of the scatterers have been obtained with MUSIC.
The original configuration of the scatterers is shown in the right image of Fig. \ref{fig:ref_for_edges}.
}
\label{fig:edges5noise}
\end{figure}

\begin{figure}[t]
\centering
\begin{tabular}{ccc}
\includegraphics[scale=0.25]{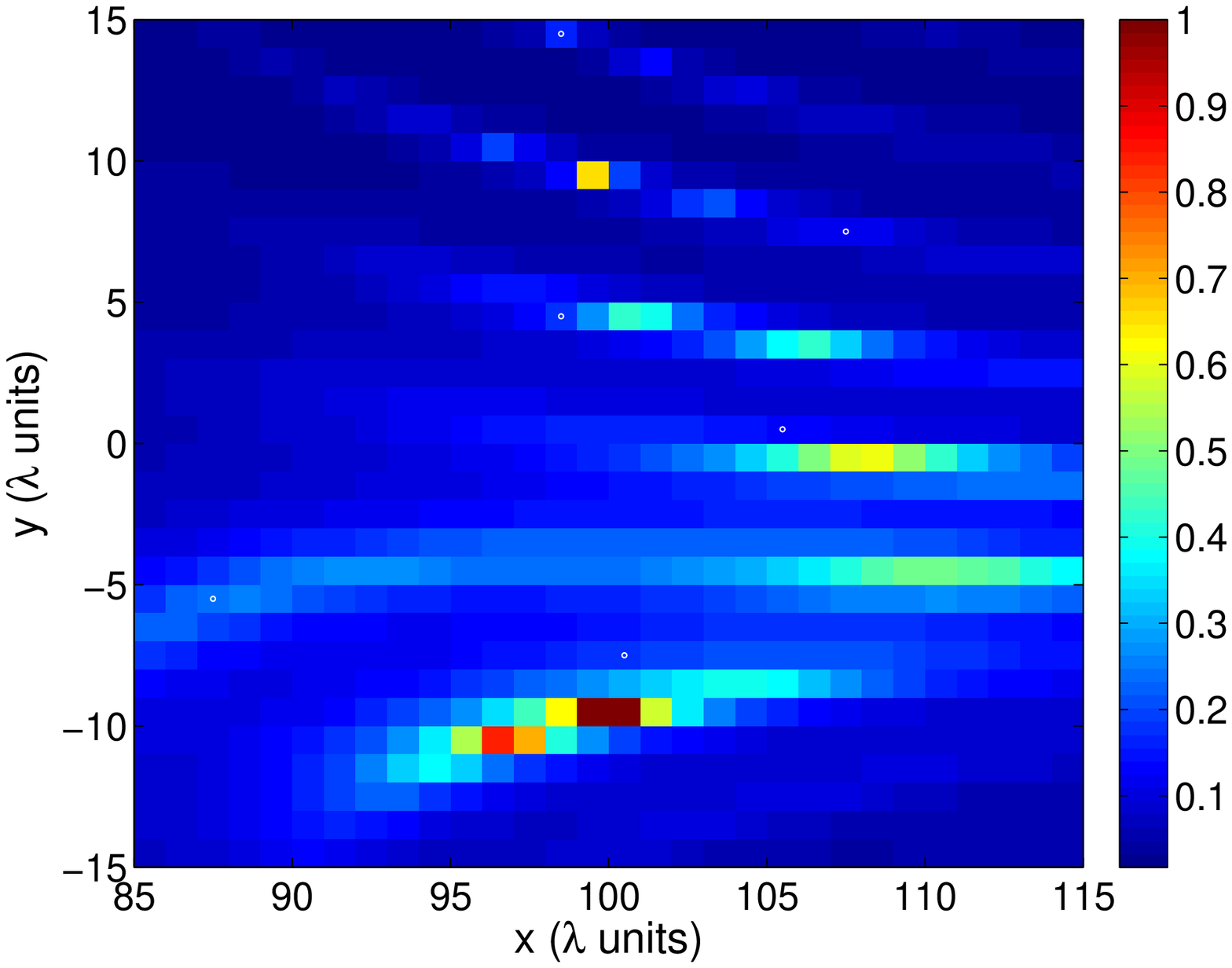} & 
\includegraphics[scale=0.25]{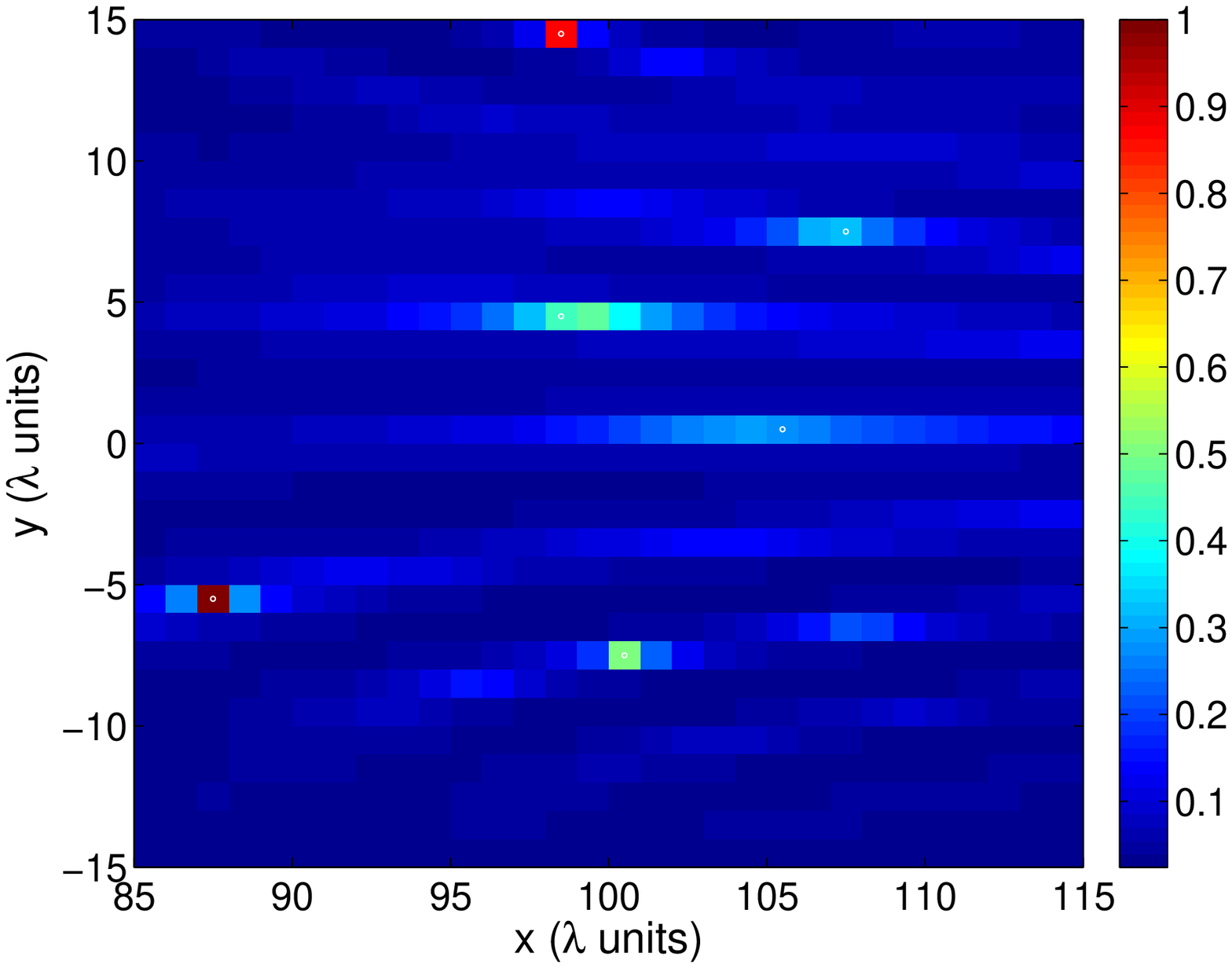} & 
\includegraphics[scale=0.25]{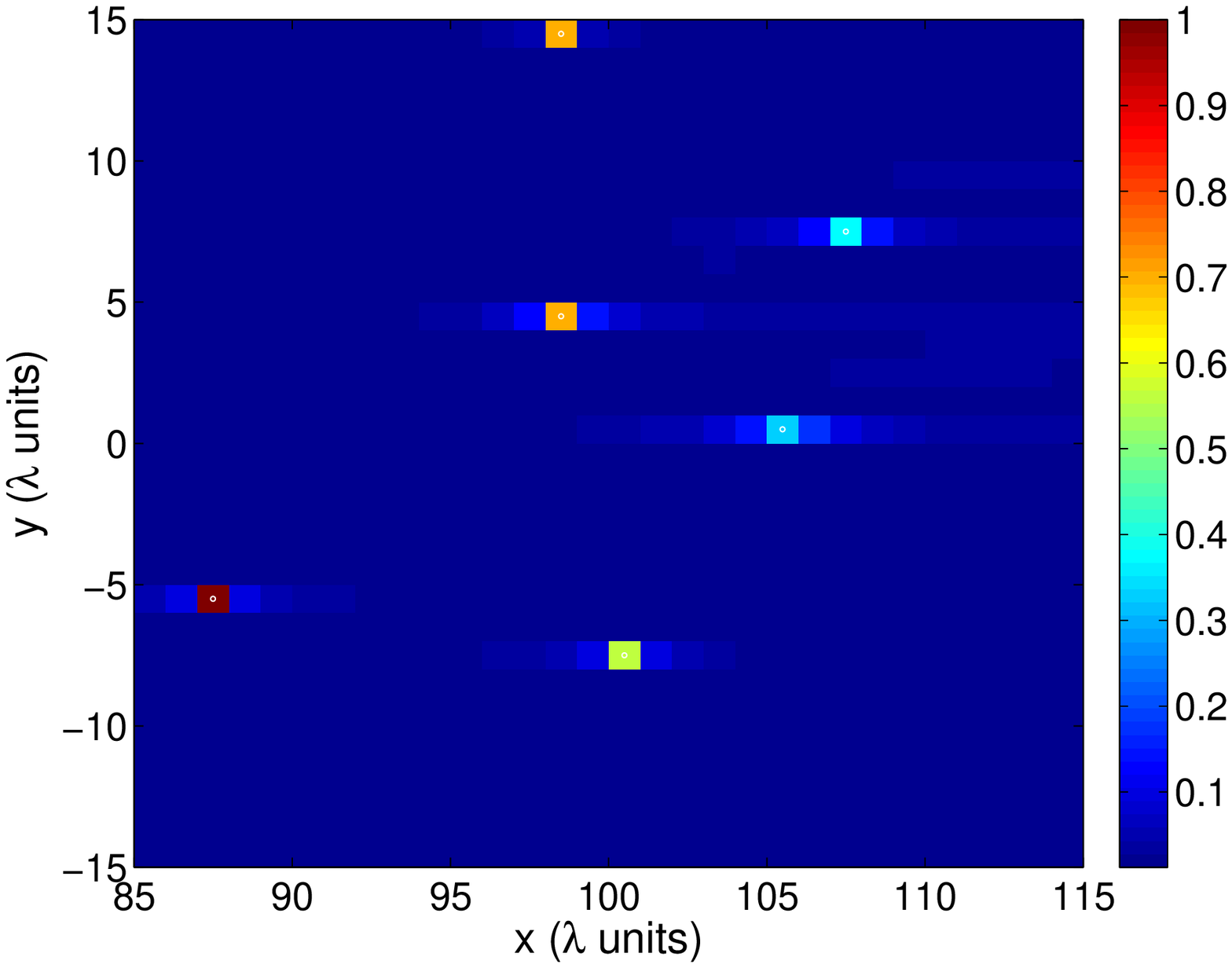} \\
\includegraphics[scale=0.25]{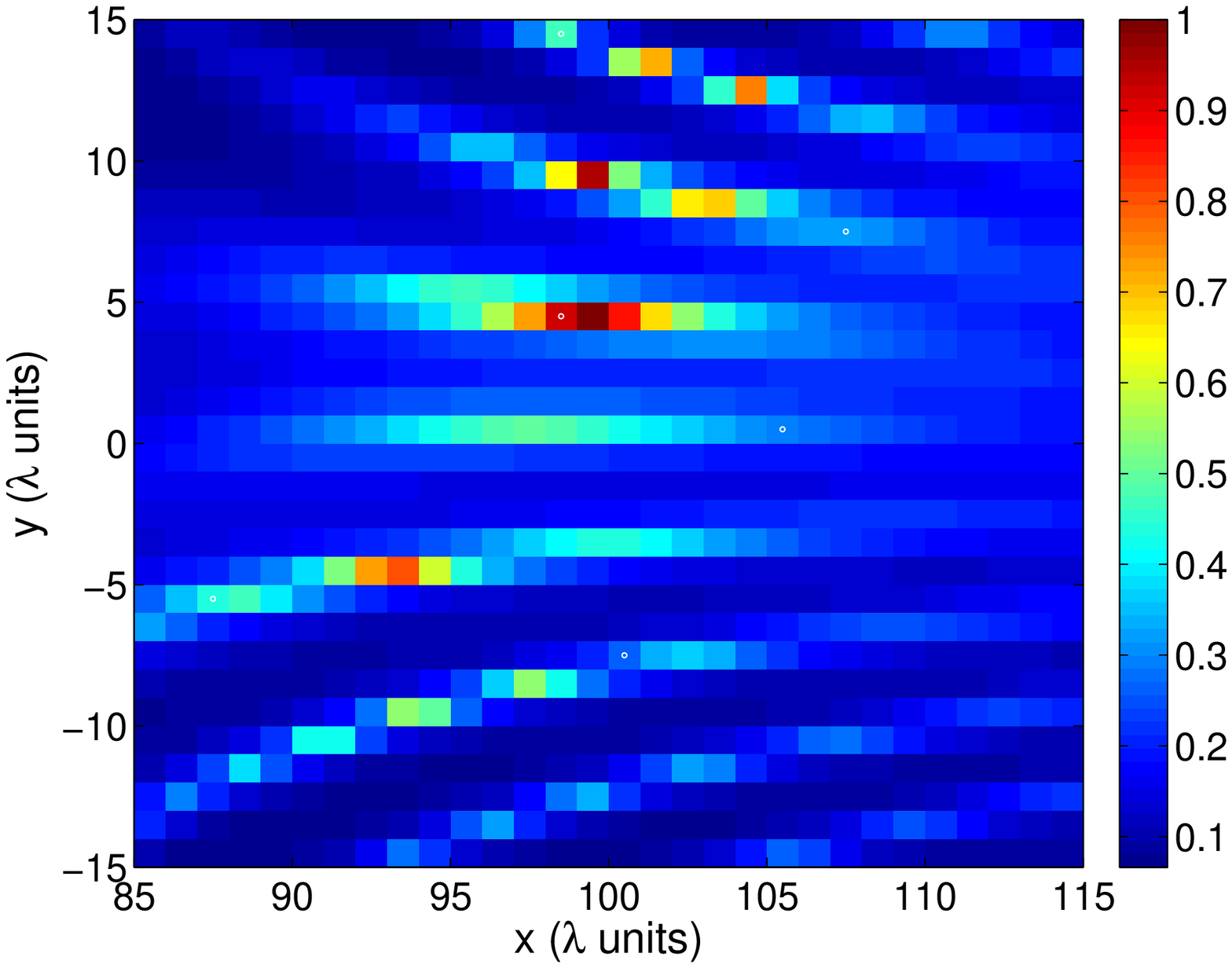} & 
\includegraphics[scale=0.25]{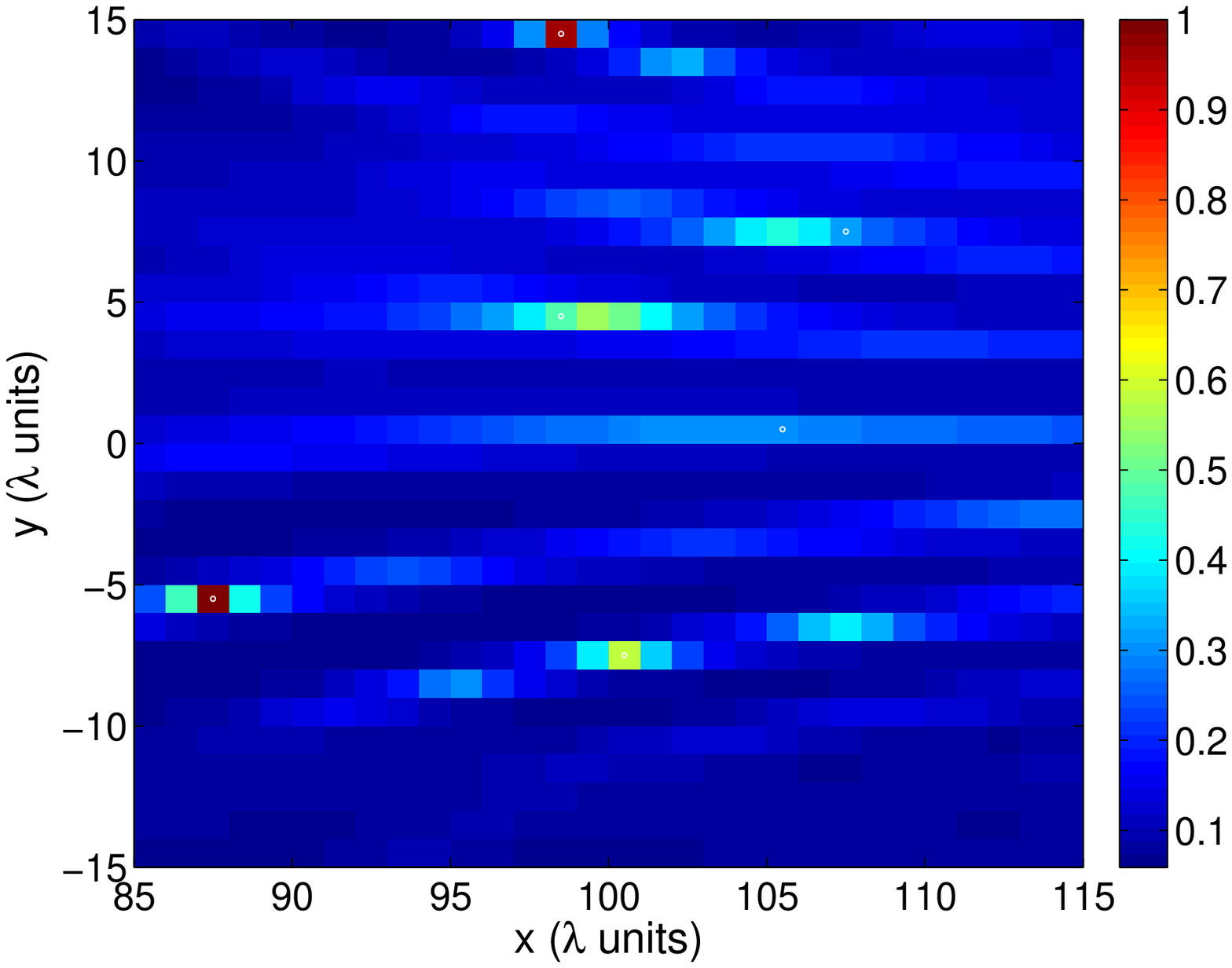} & 
\includegraphics[scale=0.25]{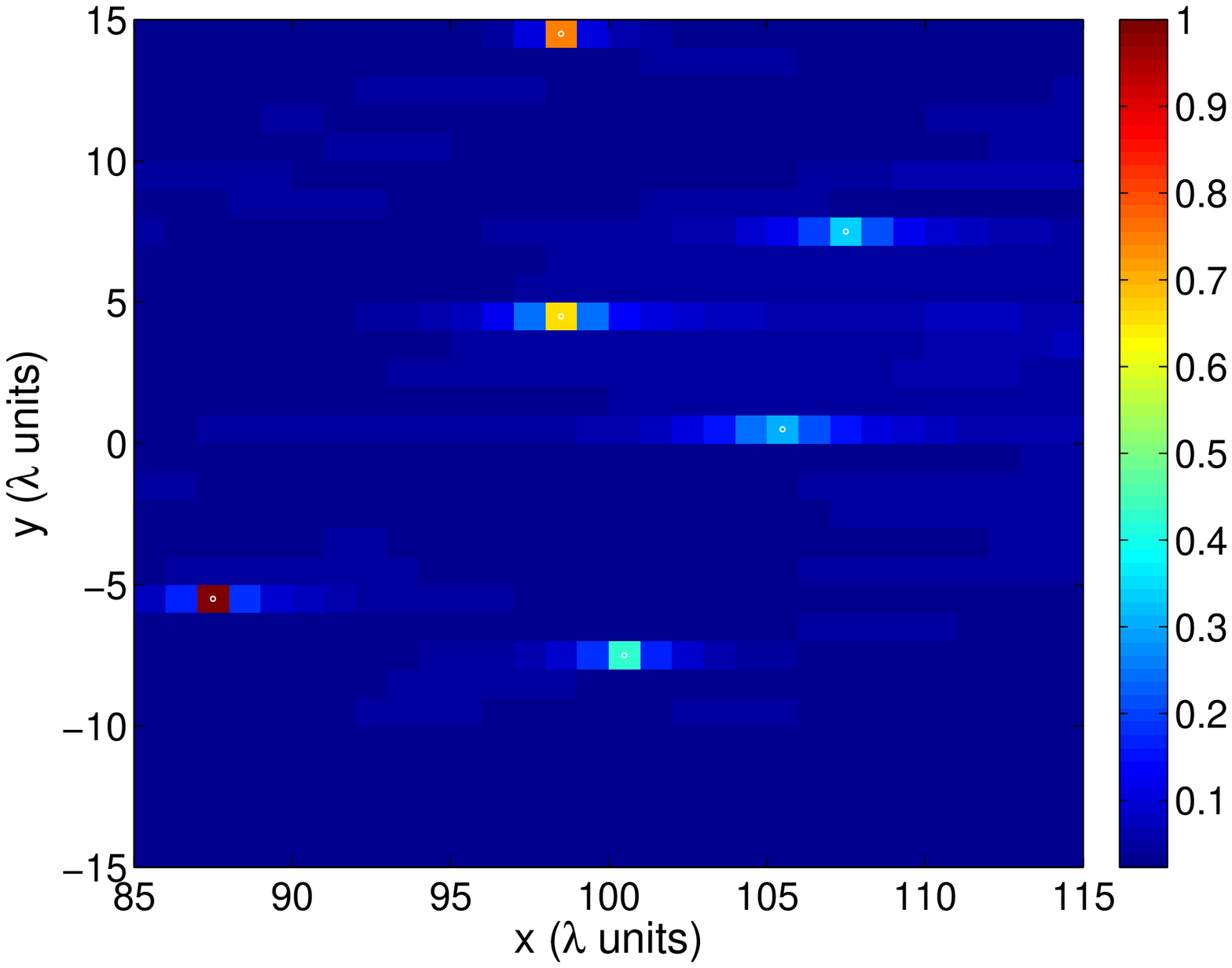}
\end{tabular}
\caption{Partial illumination from the edges of the array
with $10 \%$ (top row) and $20 \%$ (bottom row) of noise added to the data. 
$4$ (left column), $12$ (middle column) and $24$ (right column) transducers at each edge of the array illuminate the image window.
 The locations of the scatterers have been obtained with MUSIC.
The original configuration of the scatterers is shown in the right image of Fig. \ref{fig:ref_for_edges}.
}
\label{fig:edges10and20noise}
\end{figure}

To verify the robustness of the proposed illumination strategy with respect to additive noise we show in Figure
\ref{fig:edges5noise} the images obtained with MUSIC when $5 \%$ of noise is added to the data, and in Figure \ref{fig:edges10and20noise} the images obtained with MUSIC when $10 \%$ of noise (top row) and $20 \%$ of noise (bottom row) in added to the data. In Figure \ref{fig:edges5noise} we show from left to right and from top to bottom the images obtained using $4$, $8$, $12$, $16$, $20$ and $24$ transducers at each edge of the array. We see that
$16$ transducers at each edge of the array are enough to locate the scatterers accurately when $5 \%$ of noise is added to the data. In Figure \ref{fig:edges10and20noise} we see, as expected, that the higher the noise, the more transducers we need to obtain good images. The left, middle and right columns show the images obtained with
$4$, $12$ and $24$ transducers at each edge of the array, respectively.

\section{Conclusions}
\label{sec:conclusions}

We give a novel approach to imaging localized scatterers from intensity-only measurements.
The proposed approach relies on the evaluation of the {\em time reversal matrix} which, we
show, can be obtained from the total power recorded at the array using an appropriate illumination
strategy and the polarization identity. Once the {\em time reversal matrix} is obtained,
the imaging problem can be reduced to one in which the phases are known and, therefore,
one can use phase-sensitive imaging methods to form the images. These methods are very efficient, 
do not need prior information about the desired image, and guarantee the exact solution in the
noise-free case. Furthermore, they are robust with respect to noise. 

At the algorithmic level, a key property of the proposed approach is that it significantly reduces the 
computational complexity and storage consumption compared to convex approaches that replace the 
original vector problem by a matrix one \cite{CMP11,Candes13} and, therefore, create optimization 
problems of  enormous sizes. With our approach, the algorithms keep the original $K$ unknowns of 
the imaging problem, where $K$ is the number of pixels of the sought image, and hence, images of 
larger sizes can be formed.

As recording all the intensities that are needed for obtaining the {\em time reversal matrix} can be cumbersome,
we also give two solutions that simplify the data acquisition process. They greatly reduce the number of illuminations
needed for the proposed imaging strategy, but they increase the sensitivity to noise.
We illustrated the performance of the proposed strategy with various numerical
examples.


\end{document}